\documentclass[preprint,11pt,flenq,3p]{elsarticle}

\usepackage{stmaryrd} 
\usepackage{lineno}

\usepackage{hyperref}
\hypersetup{
    colorlinks=true,
    linkcolor=blue,
    filecolor=magenta,      
    urlcolor=cyan,
}
\usepackage{setspace}
\usepackage{ragged2e}
\usepackage{geometry}
\usepackage{multicol}
\usepackage{balance}       
\usepackage{graphicx}      
\usepackage{svg}
\usepackage{xcolor}
\usepackage{changepage}
\usepackage{caption}    
\usepackage{subfig} 
\usepackage[export]{adjustbox} 
\usepackage{multicol}
\usepackage{bm}
\usepackage{physics}
\usepackage{amsmath, array}
\usepackage{amssymb, mathtools}
\usepackage{verbatim} 
\usepackage{rotating}
\usepackage{indentfirst}
\usepackage{nicematrix}
\usepackage{caption}
\usepackage{multirow}
\usepackage{floatrow}
\usepackage{booktabs}   
\usepackage{enumitem}
\usepackage[nameinlink]{cleveref}
\usepackage{titletoc}
\hypersetup{urlcolor=blue, colorlinks=true, citecolor = blue}  
\hypersetup{pdfauthor={Name}}
\usepackage{cleveref}
\bibpunct{\textcolor{blue}{[}}{\textcolor{blue}{]}}{,}{a}{,}{,}
\usepackage{siunitx}
\usepackage{etoolbox}
\usepackage{xspace}
\usepackage{tikz}
\usetikzlibrary{arrows.meta, angles, quotes}

\newcommand{\pignn}{PI-GNN\xspace}          
\newcommand{\Econ}{E_{\mathrm{inc}}/E_{\mathrm{mat}}}
\newcommand{\amiya}[1]{{\textcolor{black}{#1}}}
\newcommand{\aashay}[1]{{\textcolor{black}{#1}}}
\begin{document}
\begin{frontmatter}
\title{A variational physics-informed graph neural network for heterogeneous solid mechanics} 

\author[inst1]{Aashay Rajan Yadav}
\affiliation[inst1]{organization={Mechanics of Materials Lab, Department of                  Mechanical Engineering},
            addressline={Indian Institute of Technology Madras}, 
            city={Chennai},
            postcode={600036}, 
            state={Tamil Nadu},
            country={India}
            }
\author[inst1]{Amiya Prakash Das\corref{cor}}
\ead{apdas@connect.ust.hk}
\author[inst1]{Ratna Kumar Annabattula\corref{cor}}
\ead{ratna@iitm.ac.in}
\cortext[cor]{Corresponding author}
\begin{abstract}
Stress localization in heterogeneous solids is governed by the bimaterial interface, where the displacement field remains $C^0$-continuous, while in-plane stresses jump due to the stiffness mismatch. Coordinate-based physics-informed neural networks (PINNs) represent this jump via a prescribed regularization width or a weighted interface penalty, making their accuracy sensitive to how phase-contrast changes are handled. This work presents a variational, label-free physics-informed graph neural network (\pignn) in which the heterogeneity is carried by the discretization rather than by the trial field. The solver operates on a conforming adaptive mesh graph, assigns constitutive behavior per element, and minimizes the discrete total potential energy as a single unweighted objective in which only first derivatives appear. The discrete energy on piecewise-linear elements coincides with the finite element (FE) Ritz functional. Dirichlet conditions are enforced by construction, with no penalty term, no interface weight, and no prescribed transition width. Using one fixed architecture, optimizer, and loss across small-strain elasticity and finite-strain Neo-Hookean hyperelasticity in two and three dimensions, the von Mises error remains below $3.58\%$ across a stiffness-contrast sweep spanning $(E_{\mathrm{inc}}/E_{\mathrm{mat}}\in[10^{-2},10^{2}])$, where a strong-form PINN degrades to $5.58\%$, and its displacement error reaches $7.66\%$ against $0.49\%$ for the \pignn. A trained network halves the ($\sigma_{xx}$) error of an energy-based PINN (5.01\% versus 10.94\%). Training cost exceeds a single FE solve by more than an order of magnitude, so the construction is a variationally consistent, penalty-free interface representation for parametric surrogates and inverse identification rather than a replacement for a one-off FE analysis.
\end{abstract}

\end{frontmatter}

\section{Introduction}\label{intro}
\noindent \amiya{Material heterogeneity arising from inclusions, voids, and bonded interfaces governs stress localization and failure nucleation, so accurate prediction of localized fields is central to material design. The representational difficulty lies at the bimaterial interface: the displacement field is $C^0$-continuous, and traction is continuous, whereas in-plane stresses jump across the stiffness mismatch. Smoothing that jump misrepresents the stress gradients that govern damage nucleation~\citep{henkes2022physics}. The finite element method (FEM) addresses this structurally by assigning material properties to elements on a conforming mesh.}

\amiya{Data-driven surrogates address this problem space for constitutive learning~\citep{Kirchdoerfer2016}, stress recovery~\citep{Go2025}, and composite homogenization~\citep{Maurizi2022,Xia2025}. Graph neural networks (GNNs) are mesh-native: the FE mesh becomes a graph whose message-passing layers propagate information along nodal connectivity~\citep{SanchezGonzalez2020,Zhao2024}, enabling learnable physics engines for deformation and crack propagation~\citep{Zhou2024,Wang2025}. Supervised training, however, requires large FE datasets, so the labeling burden is reassigned rather than removed. Physics-informed neural networks (PINNs) eliminate labels by enforcing governing equations directly~\citep{Raissi2019,Haghighat2021,Cuomo2022,henkes2022physics,Hu2024,ren2024mixed}, but carry two costs: a multi-objective loss requiring \emph{ad hoc} balancing and noisy second derivatives at steep stress gradients~\citep{Bai2023}. The deep energy method (DEM) removes both by minimizing the total potential energy as a single scalar objective with only first derivatives~\citep{Samaniego2020,nguyen2020deep,FuhgBouklas2022,huang2024geometrically}.}

\amiya{Neither the residual nor the energy formulation addresses the representational bias of the coordinate-based MLP (Multi-Layer Perceptron) used as the trial field. A single smooth network cannot reproduce a stress jump, so the contrast is regularized through a prescribed transition width~\citep{henkes2022physics}. Domain decomposition assigns separate networks per phase, but compatibility and traction matching re-enter as penalty terms~\citep{jagtap2020extended,jagtap2020conservative,sarma2024interface}. Alternative trial spaces\textemdash hp-VPINNs~\citep{kharazmi2021hp}, mixed-form PINNs~\citep{Rezaei2022,ren2024mixed}\textemdash still interpolate from a global coordinate map with no information about where phases meet. Evaluating the loss on an FE discretization~\citep{Zhang2025,Xiong2025,Wu2026,Rezaei2025} lets the mesh carry geometry and material assignment, but the trial field remains a dense coordinate-to-value map in which neighboring nodes exchange nothing during the forward pass.}

\amiya{GNNs close this gap by moving the discretization into the architecture. When the trial field is parametrized on the mesh graph, each message-passing step updates a nodal state based on its neighbors. Hence, a stiffness contrast between adjacent elements is visible during the forward pass, not only in the assembled objective. Existing energy-based graph solvers remain single-phase:~\citet{gao2022physics} and~\citet{he2023use} minimize variational or energy functionals on mesh graphs under small and finite strain, respectively;~\citet{dalton2023physics} scales to three-dimensional (3D) geometry with anisotropic hyperelasticity; and further applications address thermal simulation~\citep{Wurth2024}, phase-field fracture~\citep{Feng2025}, laminated shells~\citep{Hu2026}, and elastohydrodynamic lubrication~\citep{Brumand-Poor2025}. The most direct two-phase graph work~\citep{Garban2025} trains on FE stress labels with equilibrium as a weighted penalty. Label-free two-phase solvers remain coordinate-based~\citep{henkes2022physics,sarma2024interface,Rezaei2022}, representing the interface through a penalty or transition width. What has not been reported is a variational graph solver, trained without labels, in which material heterogeneity is carried by the discretization itself\textemdash with no interface penalty and no transition width\textemdash under a single construction covering linear elasticity and finite-strain hyperelasticity.}

\amiya{This paper presents a physics-informed graph neural network (\pignn) that operates on a conforming mesh graph and minimizes the discrete total potential energy. Constitutive behavior is assigned per element, so the bimaterial interface emerges from the stiffness contrast between adjacent elements exactly as in FEM. Because the discrete energy on $\mathcal{P}_1$ elements coincides with the FE Ritz functional, the converged FE field is the exact minimizer of the training objective; FEM is therefore the benchmark, with training and inference costs reported separately.}

\noindent The main contributions are:
\begin{itemize}
    \item \textbf{Heterogeneity carried by the discretization.} A single unweighted energy functional is minimized across phases with no interface penalty or transition width. Across a sweep $\Econ\in[10^{-2},10^{2}]$ with all settings frozen, the von Mises error stays below 3.58\% and displacement error below 0.5\%; a strong-form PINN with fixed transition width reaches 5.58\% and 7.66\% (\Cref{sec:pinn_vs_pignn_inclusion}).
    \item \aashay{\textbf{Isolation of the message-passing contribution.} At a fixed training budget, message passing lowers the error by \aashay{19-43\% across a mesh sweep}, with the margin widening as the mesh coarsens. At convergence, both variants settle on a similar loss, proving that message passing is a faster approach to the shared Ritz minimizer (\Cref{sec:message_passing}).}
    \item \textbf{Benchmarking on a shared unseen mesh.} The \pignn roughly halves the $\sigma_{xx}$ error of an energy-based PINN (5.01\% versus 10.94\%; 7.57\% versus 15.08\%), with errors measuring transfer to an unseen discretization (\Cref{sec:single_hole,sec:pinn_vs_pignn_inclusion}).
     \item \textbf{Robustness under one fixed model.} All verification and application examples\textemdash smooth and re-entrant inclusions, combined loading, and 3D torsion\textemdash use identical architecture, hyperparameters, and loss, with displacement $L^2$ errors below \aashay{1.5\% in both 2D and 3D (\Cref{sec:verification,sec:applications}).}
\end{itemize}

\amiya{The remainder of the paper is organized as: \Cref{sec:methodology} sets out the kinematics, constitutive relations, and variational GNN formulation. \Cref{sec:verification} reports verification and isolates the message-passing effect. \Cref{sec:applications} applies the framework to re-entrant inclusions, and 3D torsion; \Cref{sec:conclusion} closes with findings and future directions.}

\section{Preliminaries: Kinematics and \pignn Formulation}\label{sec:methodology}
\noindent This section presents the theoretical framework and formulation of the \pignn solver. We first describe the governing equations and constitutive relations, followed by the variational principle of minimum potential energy, which serves as the physical basis for the solver. Finally, we describe the graph-based discretization, neural network architecture, and training scheme in \Cref{sec:pignn_formulation}.

\subsection{Governing equations and material models}\label{sec:material_energy}
\noindent The body occupies a domain \(\Omega\) with boundary \(\partial\Omega = \Gamma_u \cup \Gamma_t\), along which the displacement \((\Gamma_u)\) and traction \((\Gamma_t)\) boundary conditions are prescribed, respectively. Given the displacement field \(\vb{u}\), the kinematic quantities are:
\begin{equation}
  \mathbb{F} = \mathbb{I} + \nabla\vb{u}, \quad
  J = \det\mathbb{F}, \quad
  \bm{\varepsilon} = \dfrac{1}{2}\!\left(\nabla\vb{u} + \nabla\vb{u}^{\top}\right), \quad
  \mathbb{C} = \mathbb{F}^{\top}\mathbb{F}, \quad
  I_1 = \operatorname{tr}(\mathbb{C}),
  \label{eq:kinematics}
\end{equation}
where \(\mathbb{F}\) and \(\mathbb{C}\) are the deformation gradient and right Cauchy-Green tensor, respectively. 

For a linear elastic (LE) material under the small-strain assumption, the strain-energy density \(\Psi\) reads:
\begin{equation}
  \Psi_{\mathrm{LE}}(\bm{\varepsilon}) = \dfrac{\lambda}{2}\,(\operatorname{tr}\bm{\varepsilon})^{2} + \mu\,(\bm{\varepsilon}\!:\!\bm{\varepsilon}),
  \label{eq:psi-le}
\end{equation}
whereas for a finite strain hyperelastic material, we use the Neo-Hookean (NH) material model,
\begin{equation}
  \Psi_{\mathrm{NH}}(\mathbb{F}) = \frac{\mu}{2}\,(I_1 - 3) - \mu\ln J + \frac{\lambda}{2}\,(\ln J)^{2}.
  \label{eq:psi-nh}
\end{equation}
The Lam\'e parameters related to the Young's modulus \((E)\) and Poisson's ratio \((\nu)\) as
\begin{equation}
  \mu = \dfrac{E}{2(1+\nu)}, \qquad
  \lambda = \dfrac{E\nu}{(1+\nu)(1-2\nu)},
  \label{eq:lame_params}
\end{equation}
governs the 3D and plane strain problems. The plane stress condition eliminates the out-of-plane stress. At finite strain, this is imposed through \(\partial\Psi_{\mathrm{NH}}/\partial \mathbb{F}_{33}=0\). And at small strain, the corresponding energy uses the reduced Lam\'e parameter \(\bar{\lambda}=E\nu/(1-\nu^{2})\). Both~\Cref{eq:psi-nh,eq:psi-le} are used, but a given boundary value problem (BVP) employs one constitutive family throughout. The linear elastic energy for the small-strain examples of~\Cref{sec:verification}, and the Neo-Hookean energy for both phases in the finite-strain examples of~\Cref{sec:applications}. Within each family, the material parameters are piecewise constant, \((E_\text{mat}, \nu_\text{mat})\) \(\in\) \(\Omega_\text{mat}\) and \((E_\text{inc}, \nu_\text{inc})\) \(\in\) \(\Omega_\text{inc}\). The total potential energy is given as:
\begin{equation}
  \Pi[\vb{u}] = \int_{\Omega} \Psi\,\mathrm{d}\Omega
  - \int_{\Omega} \vb{b}\cdot\vb{u}\,\mathrm{d}\Omega
  - \int_{\Gamma_t} \vb{t}\cdot\vb{u}\,\mathrm{d}\Gamma,
  \label{eq:pe}
\end{equation}
with body force \(\vb{b}\) and prescribed traction \(\vb{t}\) on the boundary \(\Gamma_t\). 

The principle of minimum potential energy states that the equilibrium displacement is the admissible field that minimizes \(\Pi\) (\Cref{eq:pe}):
\begin{equation}
  \vb{u}^{*} = \arg\min_{\vb{u}\in\mathcal{U}} \Pi[\vb{u}],
  \qquad
  \mathcal{U} = \{\vb{u} : \vb{u} = \vb{u}_D \text{ on } \Gamma_u\},
  \label{eq:minpe}
\end{equation}
since the stationarity condition \(\delta\Pi = 0\) is the weak form of \(\nabla \cdot\bm{\sigma}+\vb{b}=\vb{0}\) together with the natural boundary condition \(\bm{\sigma}\cdot\vb{n}=\vb{t}\) on \(\Gamma_t\). The DEM framework proposed by~\citet{Samaniego2020} uses~\Cref{eq:minpe} directly as a training objective. The displacement is approximated by the output \(\vb{u}_\theta\) of a neural network with parameters \(\bm{\theta}\). The functional is reduced to an ordinary function \(\Pi(\bm{\theta})=\Pi[\vb{u}_\theta]\) minimized over \(\bm{\theta}\) with no reference solution required. The essential condition \(\vb{u}_\theta=\vb{u}_D\) on \(\Gamma_u\) is not satisfied automatically by the network and is enforced either by construction or through a penalty term. In contrast, the traction condition on \(\Gamma_t\) is recovered naturally by the energy minimization.

\subsection{\pignn solver formulation}\label{sec:pignn_formulation}
\begin{figure}[ht]
    \centering
    \includegraphics[width=1.0\linewidth]{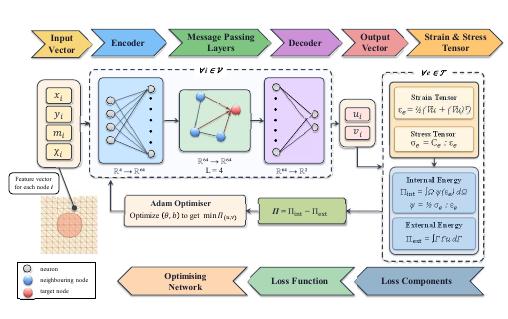}
    \caption{Overview of the \pignn solver. The mesh is a graph: each node \(i \in \mathcal{V}\) carries the feature vector \(\vb{f}_i = [\hat{X}_i,\, \hat{Y}_i,\, m_i,\, \chi_i]^{\mathsf{T}}\) (\Cref{eq:node_features}), which the encoder lifts to a 64-dimensional embedding. The \(\texttt{L = 4}\) message-passing layers (see~\Cref{fig:figure1A}) exchange embeddings along mesh edges. Element strains \(\bm{\varepsilon}_e\) and stresses \(\bm{\sigma}_e\) (linear elastic formulation shown here for simplicity) are assembled over the element set \(\mathcal{T}\) into the potential energy \(\Pi = \Pi_{\mathrm{int}} - \Pi_{\mathrm{ext}}\). The \(\Pi\) is minimized over \(\bm\theta\) using an \texttt{Adam} optimizer with no labeled data.}
    \label{fig:flowchart}
\end{figure}

\noindent The \pignn workflow is designed to solve BVPs on complex, heterogeneous domains without requiring labeled training data.~\Cref{fig:flowchart} illustrates the overall computational pipeline, which proceeds in a structured sequence: given a physical domain, we first generate an adaptive mesh and construct a graph representation \(\mathcal{G}\) (\Cref{sec:graph_construction}). We then pass spatial coordinates and material labels through a message-passing GNN to compute node-wise displacements (\Cref{sec:gnn_architecture}). Using the displacements, we compute element-wise strains and evaluate the constitutive relations to assemble the discrete potential energy \(\Pi\) (\Cref{sec:discrete_pe}). Finally, the GNN parameters \(\bm\theta\) are optimized by minimizing the energy. 

\subsubsection{Adaptive mesh generation and graph construction}\label{sec:graph_construction}
\noindent The \pignn operates on a graph \(\mathcal{G} = (\mathcal{V}, \mathcal{E})\) built from a spatially adaptive mesh of collocation points, refined where the fields are discontinuous\textemdash at the outer boundary or at material interfaces\textemdash coarsened through the quiescent bulk. Mesh generation proceeds in three stages: an adaptive size field determines the local node spacing, collocation points are placed at that spacing, and a Delaunay triangulation turns the points into a mesh graph with per-element material labels. A single geometric primitive, the signed distance \(\phi(\vb{x})\), which is negative inside the inclusion, positive in the matrix, and zero on the interface, is used. The point-in-phase tests, boundary sampling, and the size field are all derived from \(\phi\), such that the same pipeline is used for any circular (or re-entrant) inclusion and reduces to null when the inclusion interior is left unfilled.\\

\begin{figure}[h]
    \centering
    \includegraphics[width=0.75\linewidth]{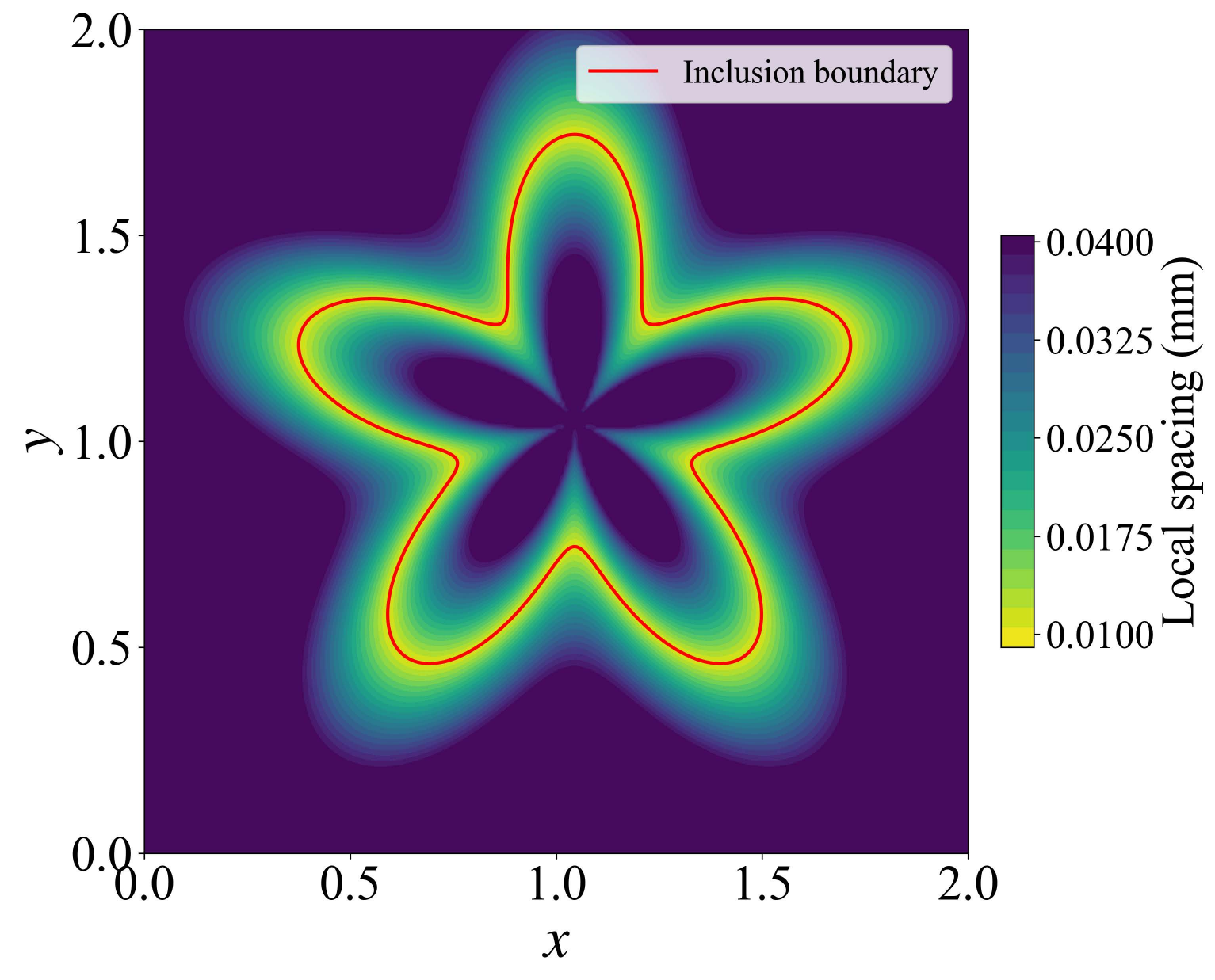}
    \caption{Adaptive collocation size field \(\rho(\vb{x})\): finest at the interface boundary (\(\rho_{\min}\)) and grading outward to the uniform far-field spacing (\(\rho_{\max}\)).}
    \label{fig:density_field_composite}
\end{figure}

\noindent \textbf{Adaptive size field:}
The field \(\rho(\vb{x})\) defines the local nodal spacing, with finer resolution prescribed geometric features and progressively coarser resolution away from them. For a feature \(k\) of characteristic length \(\ell_k\), let \(d_k = \lvert\phi_k(\vb{x})\rvert\) denote the distance from \(\vb{x}\) to its boundary. The spacing is prescribed as: 
\begin{equation}
  \rho_k(\vb{x}) = \min\!\bigl(\rho_{\min,k} + g\,d_k,\; \rho_{\max}\bigr),
  \label{eq:size_field}
\end{equation}
where \(\rho_{\min,k}\) is proportional to \(\ell_k\) specifies the minimum spacing at the feature boundary, \(g\) controls the rate of the coarsening, and \(\rho_{\max}\) denotes the far field spacing. The distance-based grading provides a gradual transition from fine to coarse resolution while avoiding an abrupt change over a prescribed refinement width. In contrast to a fixed-width refinement band, the transition width adapts to the prescribed spacing contrast and continues until \(\rho_{\max}\) is reached. The parameter \(g\) can, therefore, be selected to control the sharpness of the resolution transition or determined from the desired near-field refinement extent. When multiple features are present, the most restrictive local spacing is retained,
\begin{equation}
  \rho(\vb{x}) = \min_{k} \rho_k(\vb{x}),
  \label{eq:size_field_min}
\end{equation}
ensuring that each feature remains adequately resolved, including in regions where neighboring holes or inclusions interact. A representative spacing field is shown in~\Cref{fig:density_field_composite}.

\noindent \textbf{Collocation point assignment:}
The geometric boundaries are discretized before the inner domains are filled. Nodes are placed along each external boundary and inclusion interface according to the local spacing \(\rho\), with their positions determined by arc length. This process provides an approximately uniform distribution along smooth boundary segments while increasing the nodal density in regions of high curvatures, including re-entrant tips. The boundary and interface nodes are then fixed and treated as an exclusion region. The matrix and inclusion domains are independently filled using variable-radius Poisson-disc sampling. Pre-seeding the boundaries allows the interior sampling to automatically obtain the required clearance from all geometric interfaces. For each region, the procedure is as follows:
\begin{enumerate}
  \item \textbf{Seeding.} The fixed boundary and interface nodes are inserted into a background search grid, and the region is seeded with one accepted interior point; each inclusion is seeded separately such that none is left empty.
  \item \textbf{Candidate generation.} For a random active point, candidate points are drawn in the annulus region between one and two local spacings, with random orientation.
  \item \textbf{Acceptance.} A candidate point is accepted if it lies within the intended phase and is sufficiently separated from every existing node. The required separation is taken as the smaller of the local spacing at the candidate point and at the existing node. If this condition is satisfied for all existing nodes, the candidate point is added to the active set.
  \item \textbf{Termination.} An active point is considered retired once it yields no admissible candidates, and sampling ceases when the active set is empty, indicating that the region has reached the specified density.
\end{enumerate}

\noindent \textbf{Triangulation and material graph:}
\amiya{The complete point set, including the outer boundary, interfaces, and interiors, is triangulated using the Delaunay criterion~\citep{virtanen2020scipy}. The resulting nodes form the graph vertices \(\mathcal{V}\), while pairs of nodes sharing a triangle edge are connected by directed edges in \(\mathcal{E}\). Each edge carries a fixed geometric attribute based on the relative positions of its endpoints, \([\Delta x_{ij},\,\Delta y_{ij},\,\lVert\Delta\vb{x}_{ij}\rVert]^{\mathsf{T}}\) in 2D, with the same definition extended to \(d=3\). These attributes are computed once from the mesh and are not trained.}

\amiya{Each node is assigned a phase label; inclusion if \(\phi_k\leq0\) for any feature \(k\), including nodes on the feature boundary, and matrix otherwise. Elements are then classified according to the phases of their vertices as all-matrix, all-inclusion, or interface elements when their vertices span both labels. Because the mesh conforms to the phase interfaces, an interface element does not physically contain both materials; rather, the mixed vertex labels arise from the convention used to assign interface nodes. Such interface elements are assigned the matrix properties, consistent with the underlying conforming FE discretization and the prescribed nodal phase convention, rather than as an approximation of the material interface.}

\amiya{For a perforated plate, the same construction is used with the hole left unfilled and the remaining domain treated as a single matrix phase. The hole boundary is refined to the same field size.}

\subsubsection{Graph neural network architecture}\label{sec:gnn_architecture}
\begin{figure}[H]
    \centering
    \includegraphics[width=1.0\textwidth]{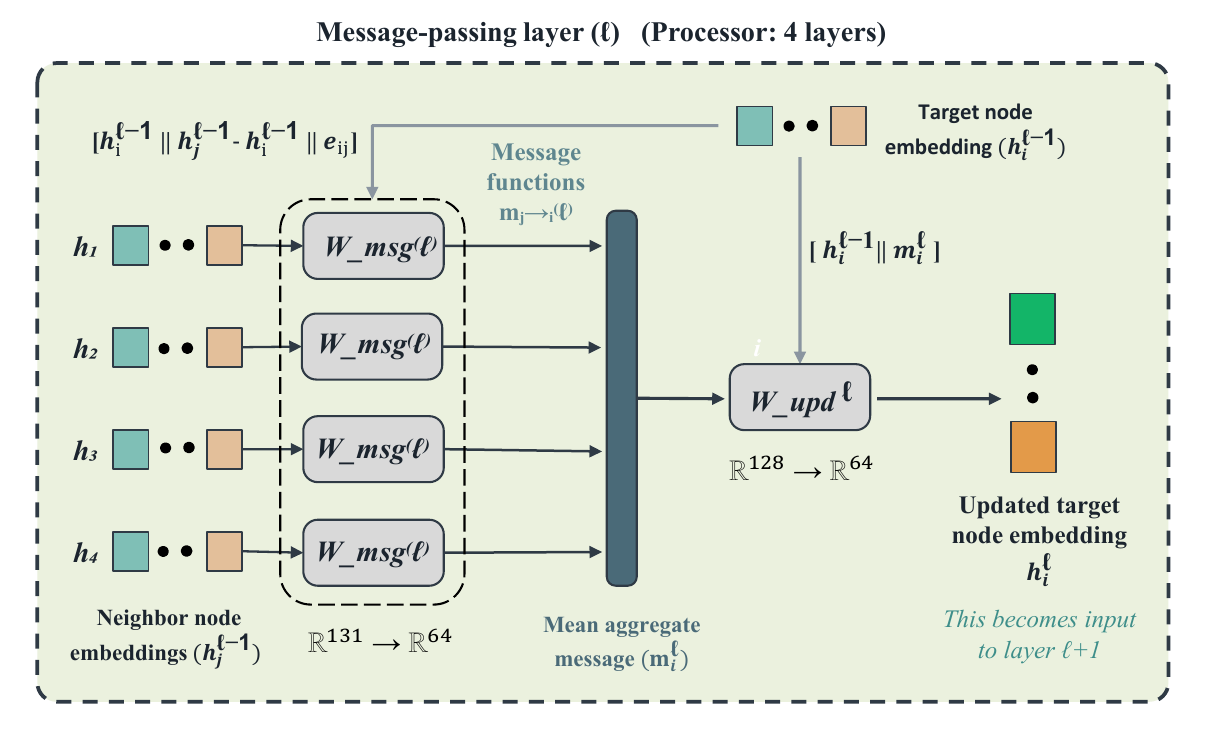}
    \caption{\aashay{Message-passing layer ($\ell$) of the processor (\(\texttt{L = 4}\) layers, \Cref{eq:message,eq:update}). On each edge ($j \to i$) the receiver embedding, the neighbor-minus-receiver difference, and the edge attribute are concatenated, $([\vb{h}_i^{(\ell-1)} \;\|\; \vb{h}_j^{(\ell-1)} - \vb{h}_i^{(\ell-1)} \;\|\; \vb{e}_{ij}] \in \mathbb{R}^{131})$, and mapped by shared weights $(\vb{W}_{\mathrm{msg}}^{(\ell)})$ to $(\vb{m}_{j \to i}^{(\ell)} \in \mathbb{R}^{64})$; the four boxes are one weight set reused edge-wise. The neighborhood mean ($\vb{m}_i^{(\ell)}$) is concatenated with ($\vb{h}_i^{(\ell-1)}$) and mapped by ($\vb{W}_{\mathrm{upd}}^{(\ell)}$) to ($\vb{h}_i^{(\ell)}$).}}
    \label{fig:figure1A}
\end{figure}

\noindent The GNN solver operates on the graph representation \(\mathcal{G}\) using an Encoder-Processor-Decoder architecture \citep{dalton2022emulation,dalton2023physics}. The encoder maps the nodal features and spatial coordinates to latent representations, which are subsequently updated through local message passing to capture neighborhood interactions. The decoder then maps the processed latent representations to nodal displacements, with the prescribed boundary conditions enforced exactly.
\begin{enumerate}
    \item \textbf{Encoder}: Each node \(i \in \mathcal{V}\) is assigned a normalized feature vector:
    \begin{equation}
      \vb{f}_i = \bigl[\hat{X}_i,\; \hat{Y}_i,\; m_i,\; \chi_i \bigr]^{\mathsf{T}} \in \mathbb{R}^4,
      \label{eq:node_features}
    \end{equation}
    where \(\hat{X}_i, \hat{Y}_i\) are normalized spatial coordinates, \(m_i \in \{0,1\}\) is the material-phase label, and \(\chi_i \in \{0,1\}\) indicates whether the node \(i\) lies on the matrix–inclusion interface. For the 2D formulation shown here; the node feature vector is \([\hat{\vb{X}}_i,\; m_i,\; \chi_i]^{\mathsf{T}} \in \mathbb{R}^{d+2}\), where \(d=2\); the same construction extends directly to 3D with ($d=3$). The decoder subsequently outputs the displacement vector \(\tilde{\vb{u}}_i \in \mathbb{R}^d\). A shared linear layer maps the node features to a 64-dimensional latent space:
    \begin{equation}
      \vb{h}_i^{(0)} = \mathrm{ReLU}\!\left(\vb{W}_{\mathrm{enc}}\,\vb{f}_i + \vb{b}_{\mathrm{enc}}\right),
      \label{eq:encoder}
    \end{equation}
    where \(\vb{W}_{\mathrm{enc}}\) and \(\vb{b}_{\mathrm{enc}}\) are trainable parameters, and \(\mathrm{ReLU}(\cdot)\) is the rectified linear unit.

    \item \textbf{Processor}: Information is propagated across the mesh using \(\texttt{L = 4}\) message-passing layers (see~\Cref{fig:figure1A}). At each layer \(\ell = 1, \dots, \texttt{L}\), node embeddings \(\vb{h}_i^{(\ell-1)}\) are updated via localized messages from incoming neighbors \(j \in \mathcal{N}(i)\):
    \begin{equation}
      \vb{m}_{j \to i}^{(\ell)} = \aashay{\mathrm{ReLU}\!\left(\vb{W}_{\mathrm{msg}}^{(\ell)} \bigl[\vb{h}_i^{(\ell-1)} \;\|\; \vb{h}_j^{(\ell-1)} - \vb{h}_i^{(\ell-1)} \;\|\; \vb{e}_{ij}\bigr] + \vb{b}_{\mathrm{msg}}^{(\ell)}\right)},
      \label{eq:message}
    \end{equation}
    \begin{equation}
      \vb{h}_i^{(\ell)} = \mathrm{ReLU}\!\left(\vb{W}_{\mathrm{upd}}^{(\ell)} \bigl[\vb{h}_i^{(\ell-1)} \;\|\; \vb{m}_i^{(\ell)}\bigr] + \vb{b}_{\mathrm{upd}}^{(\ell)}\right),
      \label{eq:update}
    \end{equation}
    \aashay{where \(\vb{m}_i^{(\ell)} = \tfrac{1}{\lvert\mathcal{N}(i)\rvert} \sum_{j \in \mathcal{N}(i)} \vb{m}_{j \to i}^{(\ell)}\) is the mean-aggregated message, \(\|\) represents vector concatenation, and \(\vb{W}^{(\ell)}\), \(\vb{b}^{(\ell)}\) are layer-specific trainable weights and biases. \amiya{The message in~\Cref{eq:message} combines three quantities: the receiver state \(\vb{h}_i^{(\ell-1)}\), the neighbor--receiver difference \(\vb{h}_j^{(\ell-1)}-\vb{h}_i^{(\ell-1)}\), and the edge geometry \(\vb{e}_{ij}\). The state difference emphasizes local contrast between neighboring nodes, which is particularly relevant near bimaterial interfaces where material properties change across adjacent elements. The edge attribute provides the relative direction and length of the connection, allowing the same message function to operate on edges with different geometries in an unstructured mesh. Because these edge features are computed directly from the mesh, the same learned message function can also be applied when the network is evaluated on an unseen discretization. \Cref{fig:figure1A} illustrates the message construction, aggregation, and node update for a single layer \(\ell\).}}

    \item \textbf{Decoder and Boundary Enforcement}: The final latent embeddings \(\vb{h}_i^{(\texttt{L})}\) are decoded to unconstrained nodal displacements via a linear layer:
    \begin{equation}
      \tilde{\vb{u}}_i = \vb{W}_{\mathrm{dec}}\,\vb{h}_i^{(\texttt{L})} + \vb{b}_{\mathrm{dec}}, \quad \tilde{\vb{u}}_i \in \mathbb{R}^2,
      \label{eq:decoder}
    \end{equation}
    where \(\vb{W}_{\mathrm{dec}}\) and \(\vb{b}_{\mathrm{dec}}\) are trainable decoder parameters. Dirichlet boundary conditions are enforced exactly by applying a multiplicative binary mask to the unconstrained predictions:
    \begin{equation}
      \vb{u}_i = \vb{q}_i \odot \tilde{\vb{u}}_i + (\vb{1} - \vb{q}_i) \odot \vb{u}_{D, i},
      \label{eq:hard_bc}
    \end{equation}
    where \(\vb{q}_i \in \{0, 1\}^2\) is the boundary mask indicating free (\(1\)) or prescribed (\(0\)) displacement components, \(\vb{u}_{D, i}\) is the prescribed Dirichlet value, and \(\odot\) denotes element-wise multiplication. This hard enforcement guarantees exact constraint satisfaction at every epoch without penalty terms in the loss.
\end{enumerate}

\subsubsection{Discrete potential energy assembly}\label{sec:discrete_pe}
\noindent Using the conforming mesh, the continuous total potential energy functional is discretized over the computational domain. With continuous piecewise-linear (\(\mathcal{P}_1\)) elements the deformation gradient \(\mathbb{F}\) and strain \(\bm{\varepsilon}\) are constant within each constant-strain triangles (CSTs) in 2D and constant-strain tetrahedra in 3D. For piecewise-constant material properties, the internal-energy density is, therefore, constant within each element, making single-point quadrature sufficient and exact for the internal energy. The domain integrals consequently reduce to algebraic sums over elements and boundary facets. Assuming no body forces, the discretized potential energy \(\Pi(\vb{u}_\theta)\) is assembled as:

\begin{equation}
  \Pi(\vb{u}_\theta) = \sum_{e=1}^{N_e} \Psi_e\,V_e
  \;-\; \sum_{f=1}^{N_f} \vb{t}_f\cdot\vb{u}^{\,f}\,S_f,
  \label{eq:pe-discrete}
\end{equation}
where \(V_e\) is the element \(e\) area in 2D (or volume in 3D), \(\Psi_e\) is the strain-energy density calculated from the predicted nodal displacements \(\vb{u}_\theta\) (using the linear elastic/hyperelastic model). \(S_f\) denotes the edge length in 2D (or facet area in 3D) of facet \(f\) on the traction boundary \(\Gamma_t\), and \(\vb{t}_f\) denotes the prescribed traction vector on that facet. For spatially varying traction, \(\vb{t}_f\) is evaluated as the facet-average traction. Similarly, \(\vb{u}^{\,f}\) denotes the predicted displacement averaged over the nodes of the facet. For a facet-wise constant traction and \(\mathcal{P}_1\) interpolation, single-point quadrature on each facet computes the external work exactly. Under the uniaxial tension condition shown in~\Cref{fig:BC_condition}, the contribution reduces to \(t\,u_x^{\,f}\,S_f\), where \(u_x^{\,f}\) is the mean longitudinal displacement of the facet nodes. \amiya{For linear elastic $\mathcal{P}_1$ elements with global displacement vector $\vb{u}$, the discrete internal strain energy evaluates to $\sum_{e=1}^{N_e} \Psi_e V_e = \tfrac{1}{2}\vb{u}^{\mathsf{T}}\vb{K}\vb{u}$, where $\vb{K}$ is the standard assembled FE stiffness matrix, and the external work evaluates to $\vb{F}_{\mathrm{ext}}^{\mathsf{T}}\vb{u}$. The discrete functional $\Pi(\vb{u}_\theta) = \tfrac{1}{2}\vb{u}^{\mathsf{T}}\vb{K}\vb{u} - \vb{F}_{\mathrm{ext}}^{\mathsf{T}}\vb{u}$ therefore coincides identically with the FE Ritz functional, and the converged FE displacement field is the exact minimizer of the discrete training loss.} The network parameters \(\bm{\theta}\) are then optimized by minimizing this discrete potential energy as the sole training objective (\Cref{eq:loss}). No PDE-residual term or boundary penalty terms are required; the essential boundary conditions are enforced hard-coded in the output layer.
\begin{equation}
  \mathcal{L}(\bm{\theta}) = \Pi(\vb{u}_\theta).
  \label{eq:loss}
\end{equation}

\section{Numerical Verification and Reference Solver}\label{sec:verification}
\noindent This section details the verification suite and benchmarks the \pignn framework against numerical and analytical references. We first describe the standard displacement-based FE reference model, which serves as the verification baseline solver (\Cref{sec:FEA_equations}). We then define the error metrics and present three verification problems (\Cref{sec:error_metrics_verification}).

\subsection{Finite element reference model}\label{sec:FEA_equations}
\begin{figure}[ht]
    \centering
    \includegraphics[width=0.95\linewidth]{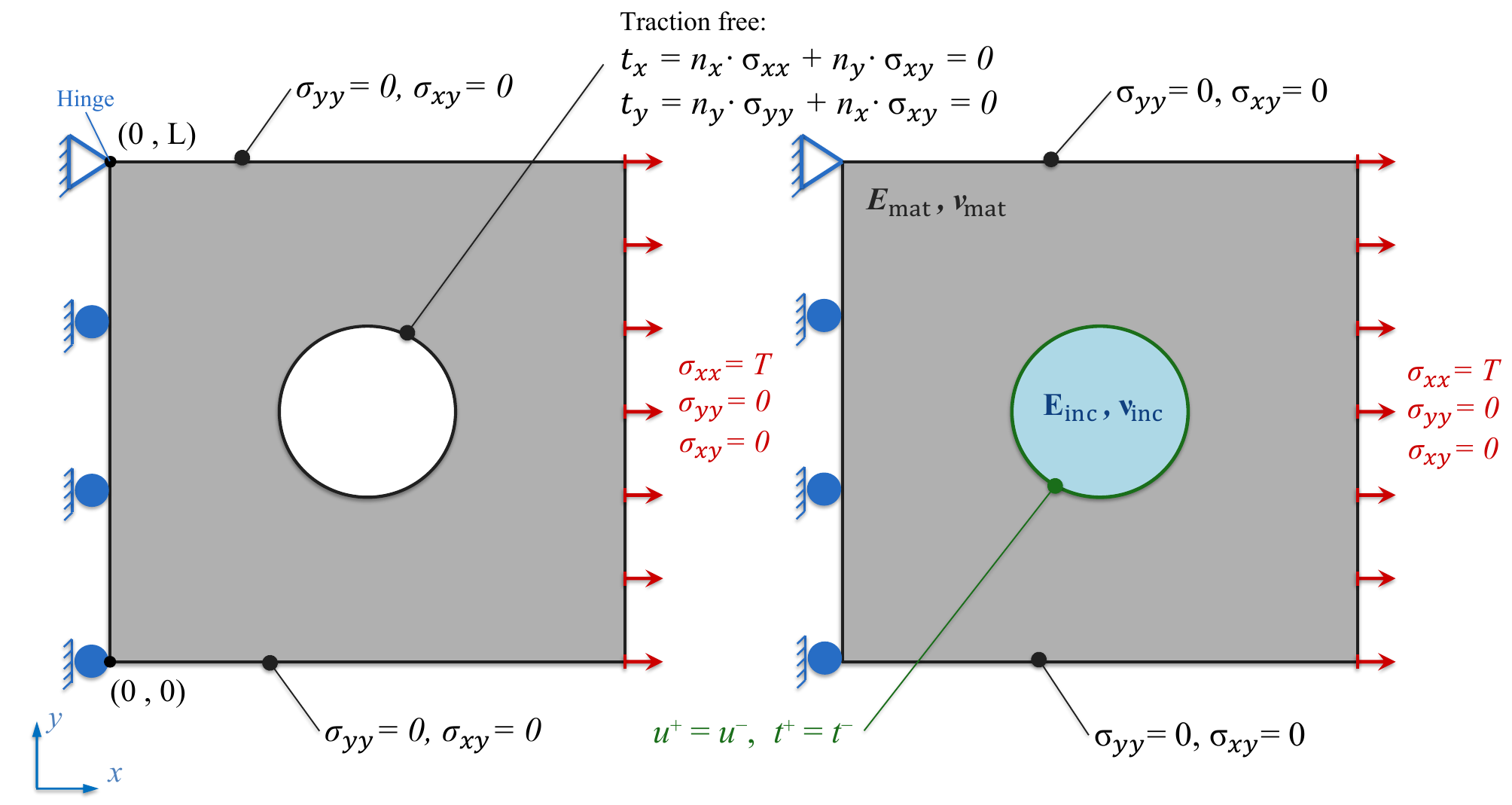}
    \caption{Boundary conditions of the representative uniaxial tension configuration, applied identically in the FE reference model and in the \pignn.}
    \label{fig:BC_condition}
\end{figure}

\noindent The predictions of the \pignn are benchmarked against a standard displacement-based FE model implemented in \texttt{FEniCSx} \citep{baratta2023dolfinx}. The solvers share the same mesh, per-phase material assignment, and boundary data, and use identical Lam\'e parameters: the shear modulus \(\mu\) throughout, with the dilatational constant taken as \(\bar{\lambda}\) for the plane stress plates and as \(\lambda\) for the 3D example. The comparison, therefore, isolates the numerical approximation introduced by the neural network from any modeling discrepancy. Converged to solver tolerance, the FE field serves as the high-fidelity benchmark throughout. Each example shares the equilibrium problem
\begin{equation}
  \nabla\!\cdot\!\bm{\sigma} = \vb{0} \quad\text{in}\ \Omega,
  \label{eq:balance_momentum}
\end{equation}
closed by the kinematics and per-phase constitutive laws and the boundary condition \(\vb{u}=\bar{\vb{u}}\) on \(\Gamma_u\) and \(\bm{\sigma}\cdot\vb{n}=\vb{t}\) on \(\Gamma_t\). The discontinuous moduli \((E,\nu)(\vb{x})\) or the Lam\'e parameters enter the element integrals directly. Because the mesh conforms to \(\Gamma_{\mathrm{int}}\), the perfect-bond interface conditions are satisfied automatically, without penalty or Lagrange-multiplier terms. \Cref{fig:BC_condition} shows the representative configuration used in the verification suite; a plate under uniaxial tension, for which the prescribed boundary conditions are:
\begin{equation}
\begin{aligned}
    u_x &= 0 && \text{on the left edge}, 
    &\qquad \vb{u} &= \vb{0} && \text{at point}\ (0,L),\\
    \bm{\sigma}\cdot\vb{n} &= (T,0)^{\!\top}
    && \text{on the right edge}.
\end{aligned}
\end{equation}
The corner constraint eliminates the remaining rigid-body mode, while all other edges are traction-free.

For the linear elastic example, the displacement is interpolated with continuous piecewise-linear (\(\mathcal{P}_1\)) Lagrange elements, and the weak form reads
\begin{equation}
  \int_{\Omega}\bm{\sigma}(\vb{u}_h):\bm{\varepsilon}(\vb{v})\,\mathrm{d}\Omega = \int_{\Gamma_t}\vb{t}\cdot\vb{v}\,\mathrm{d}\Gamma,
  \label{eq:weak_form}
\end{equation}
where \(\bm{\sigma}(\vb{u}_h)\) is the stress conjugate to \(\Psi_{\mathrm{LE}}\). The spatially varying \(\mu(\vb{x})\) and \(\bar{\lambda}(\vb{x})\) enter the bilinear form directly, requiring no modification of the formulation for homogeneous or composite domains. The resulting linear system is solved using a direct LU factorization with MUMPS \citep{amestoy2001fully}. Element stresses are subsequently recovered from \(\vb{u}_h\) using the same constitutive relation.

For the hyperelastic example, the mesh, element type, and boundary treatment are retained, while the constitutive model is replaced by the Neo-Hookean energy \(\Psi_{\mathrm{NH}}\), rendering the problem nonlinear. A total-Lagrangian formulation is adopted with respect to the reference configuration. The discrete equilibrium is obtained from the stationarity of the total potential energy, whose first variation defines the residual:
\begin{equation}
  \mathcal{R}(\vb{u}_h;\vb{v}) = \int_{\Omega}\mathbb{P}:\nabla\vb{v}\,\mathrm{d}\Omega - \int_{\Gamma_t}\vb{t}\cdot\vb{v}\,\mathrm{d}\Gamma = 0 \qquad \forall\,\vb{v}\in\mathcal{V}_0,
  \label{eq:fem_residual_nh}
\end{equation}
where \(\mathbb{P}=\partial\Psi_{\mathrm{NH}}/\partial\mathbb{F}\) is the first Piola-Kirchhoff stress and \(\mathcal{V}_0\) is the space of admissible variations vanishing on \(\Gamma_u\). The nonlinear system is solved using Newton-Raphson iteration with the consistent tangent \(\partial\mathcal{R}/\partial\vb{u}_h\) assembled by automatic differentiation of the residual. The prescribed nominal traction \(\vb{t}\) on the reference boundary is applied incrementally in equal load steps to facilitate convergence at finite strain, with the direct LU solver reused at each iteration. The Cauchy stress used for comparison is subsequently recovered by push-forward,
\begin{equation}
  \bm{\sigma} = J^{-1}\mathbb{P}\mathbb{F}^{\top}, \qquad \mathbb{P}=\mathbb{F}\mathbb{S}, \qquad \mathbb{S} = \mu\!\left(\mathbb{I}-\mathbb{C}^{-1}\right) + \lambda_{\mathrm{3D}}\,\ln J\,\mathbb{C}^{-1},
  \label{eq:cauchy_nh}
\end{equation}
where \(\mathbb{S}\) denotes the second Piola-Kirchhoff stress.

In linear elastic or hyperelastic regimes, the FE reference solution and the \pignn use the same mesh, per-element material assignment, and essential and natural boundary conditions. The comparison, therefore, isolates the discrepancy between the converged FE solution and its neural approximation without introducing discrepancies due to geometry, discretization, material assignment, or loading.

\subsection{Verification problems and error metrics}\label{sec:error_metrics_verification}
\noindent The framework is verified across two problem classes: the first investigates perforated plates with single- and three-hole configurations that introduce stress concentrations. The second investigates a material interface with a stiffness mismatch. These numerical examples probe the recovery of smooth elastic fields, sharp stress gradients, and heterogeneous response across an interface. 

Where a closed-form solution is available, it is used as the reference; otherwise, the neural predictions are compared with the converged FE solution evaluated on the \emph{same} mesh. This common discretization removes differences in mesh resolution and geometry representation from the comparison, such that the reported error isolates the neural approximation relative to the FE reference. Body forces are neglected throughout both the FE and neural network formulations, with \(\vb{b}=\vb{0}\) in~\Cref{eq:pe}.

Let \(\alpha\) denote a scalar field such as a displacement component \(u_x,u_y\), a stress component \(\sigma_{xx},\sigma_{yy},\sigma_{xy}\), or the von Mises stress \(\sigma_{\mathrm{vM}}\). Let \(\alpha_{\mathrm{FEM}}\) and \(\alpha_{\mathrm{PI\text{-}GNN}}\) denote the FE reference and \pignn predictions, respectively, evaluated at the same nodes. Accuracy is quantified using three complementary measures: the global relative \(L^2\) error, the pointwise relative error \(\text{err}_{\alpha}\), and the coefficient of determination \(R^2\). The global agreement between the two solutions is quantified using an element-area-weighted relative \(L^2\) error for each field:
\begin{equation}
    L^{2}_{\alpha} = \dfrac{\sqrt{\displaystyle\sum_{e} A_e\,\bigl(\alpha^{e}_{\mathrm{PI\text{-}GNN}} - \alpha^{e}_{\mathrm{FEM}}\bigr)^2}}{\sqrt{\displaystyle\sum_{e} A_e\,\bigl(\alpha^{e}_{\mathrm{FEM}}\bigr)^2}}\times 100\,\%,
    \label{eq:4.2}
\end{equation}
where \(A_e\) is the area of element \(e\) and \(\alpha^{e}\) the field value on that element, so that the metric weights each region by its geometric measure and is independent of mesh non-uniformity. The signed pointwise relative error at node \(i\) is evaluated as:
\begin{equation}
    \text{err}_{\alpha}(\vb{x}_i) = \dfrac{\alpha_{\mathrm{PI\text{-}GNN}}(\vb{x}_i) - \alpha_{\mathrm{FEM}}(\vb{x}_i)}{\lvert \alpha_{\mathrm{FEM}}(\vb{x}_i)\rvert}\times 100\,\%.
    \label{eq:4.1}
\end{equation}

The sign of \(\text{err}_{\alpha}\) distinguishes over-prediction (\(+\)) from under-prediction (\(-\)). To avoid the artificial amplification of the relative error in~\Cref{eq:4.1} near zero-crossings and constrained boundaries, nodes satisfying \(\lvert \alpha_{\mathrm{FEM}}(\vb{x}_i)\rvert < \delta\,\max_j \lvert \alpha_{\mathrm{FEM}}(\vb{x}_j)\rvert\) are excluded, with a relative threshold of \(\delta = 10^{-2}\). The spatial distribution of \(\text{err}_\alpha\) is summarized by its median, which indicates systematic bias toward over- or under-prediction while remaining less sensitive to isolated outliers. The coefficient of determination measures the reference field variance captured by the prediction,
\begin{equation}
    R^2_{\alpha} = 1 - \dfrac{\displaystyle\sum_{i}\bigl(\alpha_{\mathrm{FEM}}(\vb{x}_i) - \alpha_{\mathrm{PI\text{-}GNN}}(\vb{x}_i)\bigr)^2}{\displaystyle\sum_{i}\bigl(\alpha_{\mathrm{FEM}}(\vb{x}_i) - \overline{\alpha}_{\mathrm{FEM}}\bigr)^2},
    \label{eq:4.4}
\end{equation}
where \(\overline{\alpha}_{\mathrm{FEM}}\) is the mean of the reference field. For the heterogeneous examples, \(R^2_\alpha\) is additionally evaluated separately in near-field and far-field regions. A node is classified as near-field if its distance to the boundary of any inclusion is less than or equal to twice the corresponding inclusion radius; all remaining nodes are classified as far-field. This partition separates the accuracy in regions containing steep gradients near material interfaces from that in the comparatively smooth bulk.

\subsubsection{Single and multiple holes}\label{sec:single_hole}
\noindent The first verification example considers square plates with circular holes, providing a benchmark for localized stress concentrations. Three configurations of increasing geometric complexity are used to assess the \pignn. The smallest single-hole example first verifies the \pignn and the FE reference against the closed-form Kirsch solution, for which the finite plate approximates the infinite-plate idealization well (\(a=\SI{1}{\milli\meter}\), plate-to-diameter ratio of 10). The two remaining configurations, consisting of a larger single hole (\(r=\SI{4}{\milli\meter}\)) and an asymmetric three-hole plate, additionally provide a comparison with the energy-based PINN \citep{Li2021}. For neural comparisons, the solvers are trained using independent samples from the domain and subsequently evaluated, without further optimization, with a single forward pass on a common FE mesh that is not used during training. The converged FE solution on this mesh serves as the reference.

The \pignn versus PINN comparison is designed to isolate the effect of the displacement-field representation as far as possible. Both models minimize the same functional, \(\Pi=U-W_{\mathrm{ext}}\), using \texttt{Adam} with a constant learning rate of \(\eta=1\times10^{-3}\) and a fixed budget of \(30{,}000\) epochs, without problem-specific tuning. The PINN contains approximately \(21{,}000\) trainable parameters, compared with approximately \(17{,}000\) for the \pignn \texttt{L4H32}, providing comparable network capacities and avoiding an advantage to the \pignn from a larger parameter count. The \pignn represents the displacement as a discrete nodal field and evaluates the energy element-wise on the mesh. In contrast, the PINN represents a continuous coordinate-based field and evaluates the energy by Monte Carlo integration over collocation points \citep{Li2021}. Thus, both methods are compared under the same variational objective and optimization budget, while their spatial representations and energy discretizations remain distinct.\\

\noindent\textbf{Analytical benchmark:} The plate occupies \(\Omega = [-l/2, l/2]^2 \setminus \mathcal{B}_a(\vb{0})\) with side length \(l = \SI{20}{\milli\meter}\) and hole radius \(a = \SI{1}{\milli\meter}\). The resulting plate-to-diameter ratio, \(l/2a = 10\), provides a finite approximation to the infinite-plate assumption underlying the Kirsch solution. The material is linear elastic steel (\(E = \SI{210}{\giga\pascal}\), \(\nu = 0.3\)) under plane stress. A uniform tensile traction \(T = \SI{1}{\mega\pascal}\) is applied on the right edge, with the Dirichlet conditions on the left edge (see~\Cref{fig:BC_condition}). The near-field region is defined as the annulus extending from the hole boundary to an outer radius \(2a\).

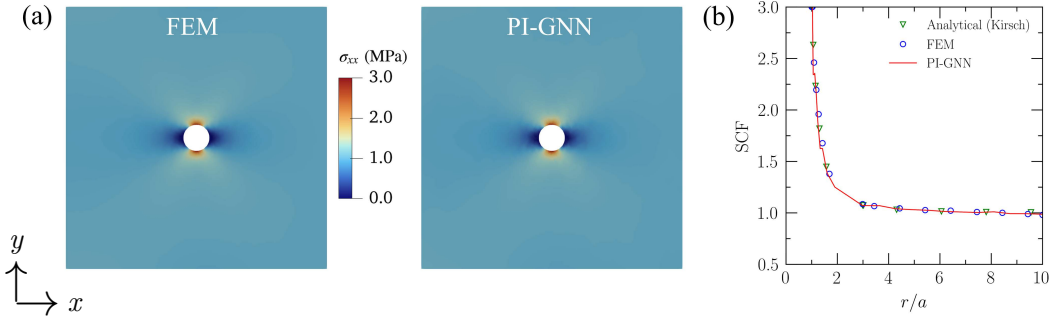
\begin{figure}[H]
    \centering
    \includegraphics[width=0.9\linewidth]{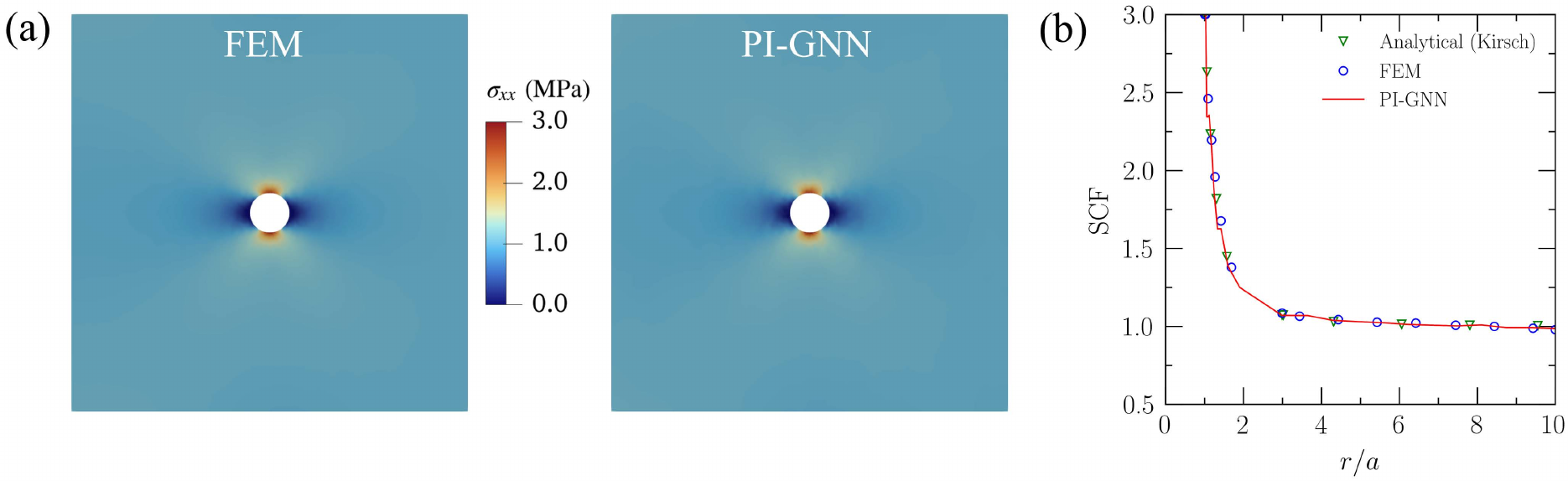}
    \put(-390,0){\LARGE\(\rightarrow\)}
    \put(-393.5,8.0){\LARGE\(\uparrow\)}
    \put(-370,1.5){\(\textit{x}\)}
    \put(-392,25){\(\textit{y}\)}
    \caption{Thin plate with a single hole; \(a = \SI{1}{\milli\meter}\), plate-to-diameter ratio 10. (a) Longitudinal stress \(\sigma_{xx}\) from the FE reference solution (left) and the \pignn prediction (right). (b) \(\sigma_{xx}\) along the transverse line of symmetry \(x = 0\) from the \pignn, the FE reference, and the analytical Kirsch solution, showing the decay of the stress-concentration factor from 3 at the hole boundary to unity in the far field.}
    \label{fig:single_hole_kirsch}
\end{figure}

\Cref{fig:single_hole_kirsch}a compares the longitudinal stress \(\sigma_{xx}\) of the FE reference solution (left) and the \pignn prediction (right). The \pignn reproduces the characteristic Kirsch distribution; \(\sigma_{xx}\) reaches \(\approx 3T\) at the hole crown (\(\theta = \pm\pi/2\)). At the horizontal margins (\(\theta = 0, \pi\)), the pressure becomes compressive with a value \(\approx -T\), and approaches the far-field value \(T\) away from the hole. \Cref{fig:single_hole_kirsch}b traces \(\sigma_{xx}\) along the transverse line of symmetry \(x = 0\). The \pignn closely follows both the FE reference and the analytical Kirsch solution, recovering the maximum stress-concentration factor of 3 at the hole boundary and its decay toward the far-field value. The field-wise metrics in \Cref{tab:r2_scores_1hole} shows good agreement, with the load-carrying fields \(u_x\) and \(\sigma_{xx}\) achieving \(R^2 = 0.99\) both globally and in the near-field. \aashay{The mesh convergence study for this case is shown in \Cref{app:mesh_conv_1_hole}.}

\begin{table}[H]
    \centering
    \small
    \setlength{\tabcolsep}{12pt} 
    
    \caption{Comparison of the field-wise error for the thin plate with a single-hole (\(a = \SI{1}{\milli\meter}\)): \pignn against the FE reference.}
    \label{tab:r2_scores_1hole}
    
    \begin{tabular}{l S[table-format=1.4] S[table-format=1.4] S[table-format=1.4]}
        \toprule \toprule
        & & \multicolumn{2}{c}{\(R^2\)} \\
        \cmidrule(l){3-4}
        Field & {\(L^{2}_\alpha\) (\%)} & {Global} & {Near-field} \\
        \midrule
        \(u_x\)         & 0.5700 & 0.9998 & 0.9992 \\
        \(u_y\)         & 0.5700 & 0.9996 & 0.9989 \\
        \(\sigma_{xx}\) & 0.9397 & 0.9963 & 0.9964 \\
        \(\sigma_\text{vM}\) & 0.8879 & 0.9963 & 0.9964 \\
        \bottomrule \bottomrule
    \end{tabular}
\end{table}

\noindent \textbf{Single hole (\(r = \SI{4}{\milli\meter}\)):} The problem is posed on the full plate (side \SI{20}{\milli\meter}), containing a central hole of radius \(r = \SI{4}{\milli\meter}\) and subjected to the roller-and-pin conditions shown in~\Cref{fig:BC_condition}. This differs from the quarter-domain formulation with symmetry conditions used in the original benchmark~\citep{Li2021}. The \pignn is trained on its own conforming mesh, refined around the hole, containing 4368 nodes and 8221 elements. In contrast, the PINN is trained using \(10{,}000\) collocation points, which are resampled every 500 epochs via near-field densification. Both are then evaluated on a common FE mesh with 1865 nodes and 3420 elements, which is not used during training. \aashay{The \pignn is evaluated by a forward pass on the mesh graph, while the PINN is evaluated by querying its coordinate network at the mesh nodes.}

\begin{figure}[H]
    \centering
    \includegraphics[width=0.9\linewidth]{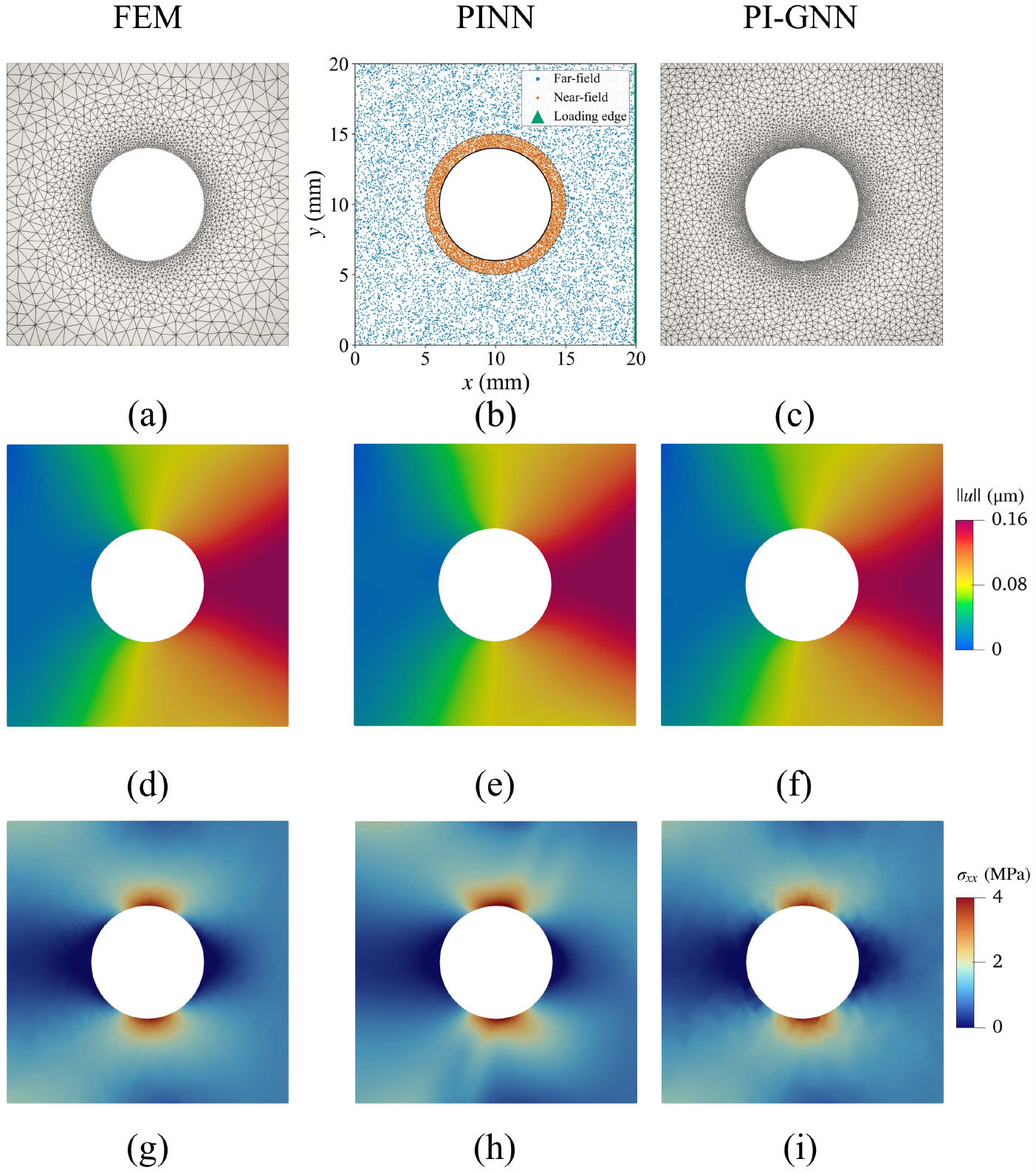}
    \caption{Discretizations and field predictions for the single-hole plate ($r = \SI{4}{\milli\meter}$). 
    \textbf{Top row:} (a) Shared FE evaluation mesh ($1865$ nodes, $3420$ elements). (b) Collocation points for the PINN energy estimate, sampled with near-field densification ($10{,}000$ points, resampled every $500$ epochs)~\citep{Li2021}. (c) \pignn training mesh ($4368$ nodes, $8221$ elements). The neural solvers are trained on (b) and (c), then evaluated on (a). 
    \textbf{Middle row:} Displacement magnitude ($\Vert\vb{u}\Vert$) fields evaluated on the FE reference mesh. (d) FE reference solution, (e) energy-based PINN prediction~\citep[re-implemented from][]{Li2021}, and (f) \pignn prediction. 
    \textbf{Bottom row:} Longitudinal stress ($\sigma_{xx}$) fields evaluated on the FE reference mesh. (g) FE reference solution, (h) energy-based PINN prediction, and (i) \pignn prediction.}
    \label{fig:pinn_vs_pignn_sigma_xx}
\end{figure}

The displacement magnitude \((||\vb{u}||)\) and longitudinal stress \((\sigma_{xx})\) fields exhibit good qualitative agreement, with all three solvers accurately resolving the characteristic crown stress concentration of approximately \(3T\) and its decay into the surrounding plate (see \Cref{fig:pinn_vs_pignn_sigma_xx}). Quantitatively, the \pignn provides better agreement with the FE reference than the PINN, achieving a global \(R^2 = 0.99\) compared with 0.96 for the PINN. The same trend is observed in the near-field, where the \pignn and PINN achieve \(R^2=0.99\) and 0.97, respectively, indicating that the additional near-field sampling of the PINN does not substantially improve its accuracy in the stress-concentration region (see \Cref{tab:pinn_vs_pignn_scores}). The relative-error distributions further distinguish the two models. The \pignn errors are narrowly distributed around zero, with a median of \(-0.41\%\) and a standard deviation of \(5.24\%\), whereas the PINN exhibits a positive bias, with a median of \(+5.07\%\), and a broader distribution with a standard deviation of \(9.69\%\) (\Cref{fig:disp_err_density}c).\\

\noindent \textbf{Three-hole configuration:} To test the resolution of interacting stress fields, an asymmetric three-hole configuration is considered in which closely spaced holes of different radii generate overlapping stress concentrations. The square plate (side \(l = \SI{20}{\milli\meter}\)) contains three circular holes of radii \(r_1 = \SI{1}{\milli\meter}\), \(r_2 = \SI{2}{\milli\meter}\), and \(r_3 = \SI{3}{\milli\meter}\) centered at \((8, 14)\), \((12, 12)\), and \((8, 7)\), respectively. The material, loading, and boundary conditions are identical to the single-hole example. The \pignn is trained on its own conforming mesh refined around the three holes, containing 4455 nodes and 8394 elements, while the PINN is trained using \(25{,}000\) collocation points resampled every 500 epochs. Both models are subsequently evaluated on a common unseen FE mesh containing 3158 nodes and 5837 elements, with the converged FE solution on this mesh serving as the reference.

\begin{figure}[H]
    \centering
    \includegraphics[width=0.9\linewidth]{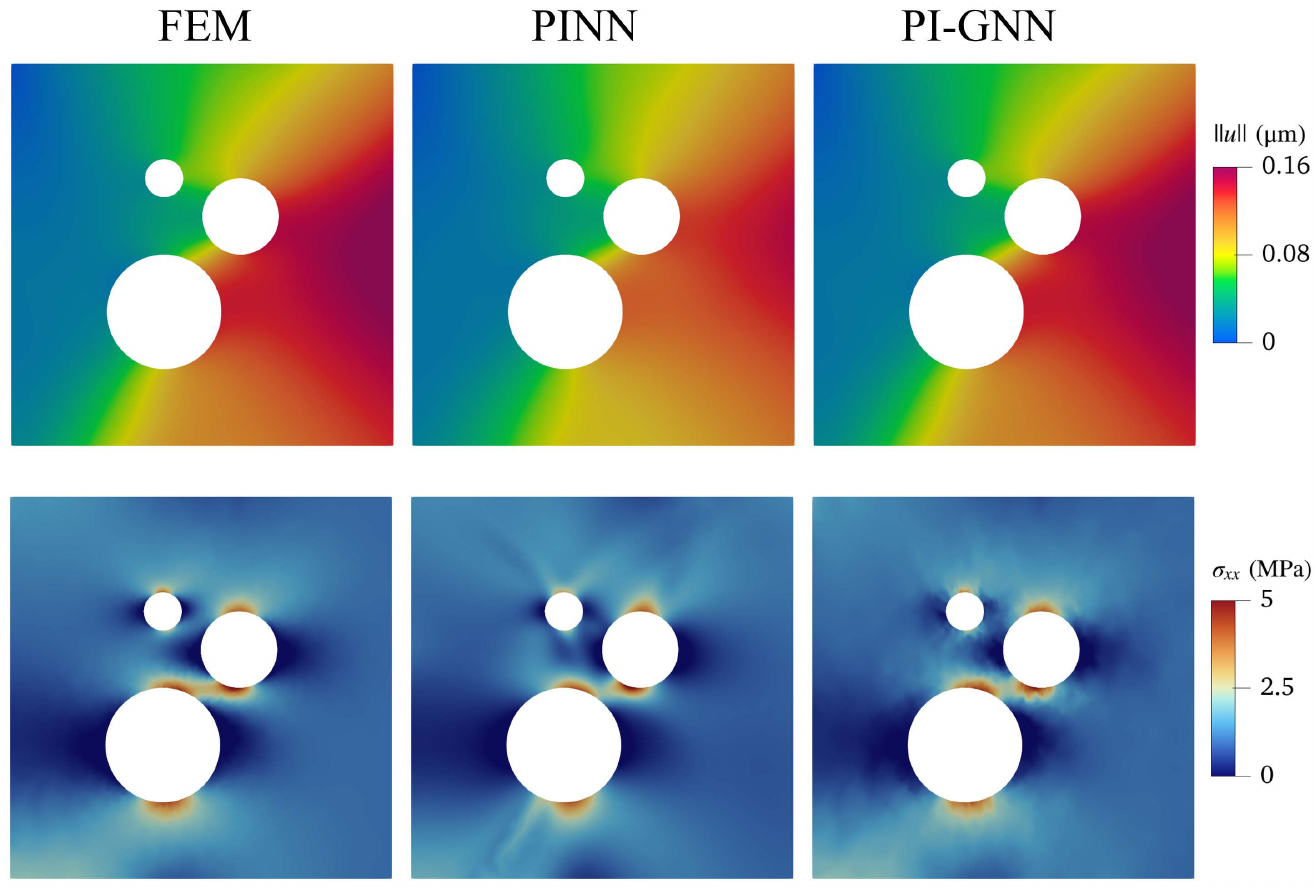}
    \caption{\textbf{Top row:} Displacement magnitude ($\Vert\vb{u}\Vert$). \textbf{Bottom row:} longitudinal stress \(\sigma_{xx}\) for the asymmetric three-hole configuration~\citep{Li2021}. The peak \(\sigma_{xx} \approx 5T\) develops in the ligament above the largest hole. (a) FE reference solution (\texttt{FEniCSx}). (b) Energy-based PINN prediction,~\citep[re-implemented here][]{Li2021}. (c) \pignn prediction.}
    \label{fig:reference_multihole_configuration}
\end{figure}

\begin{figure}[H]
    \centering
    \includegraphics[width=0.9\linewidth]{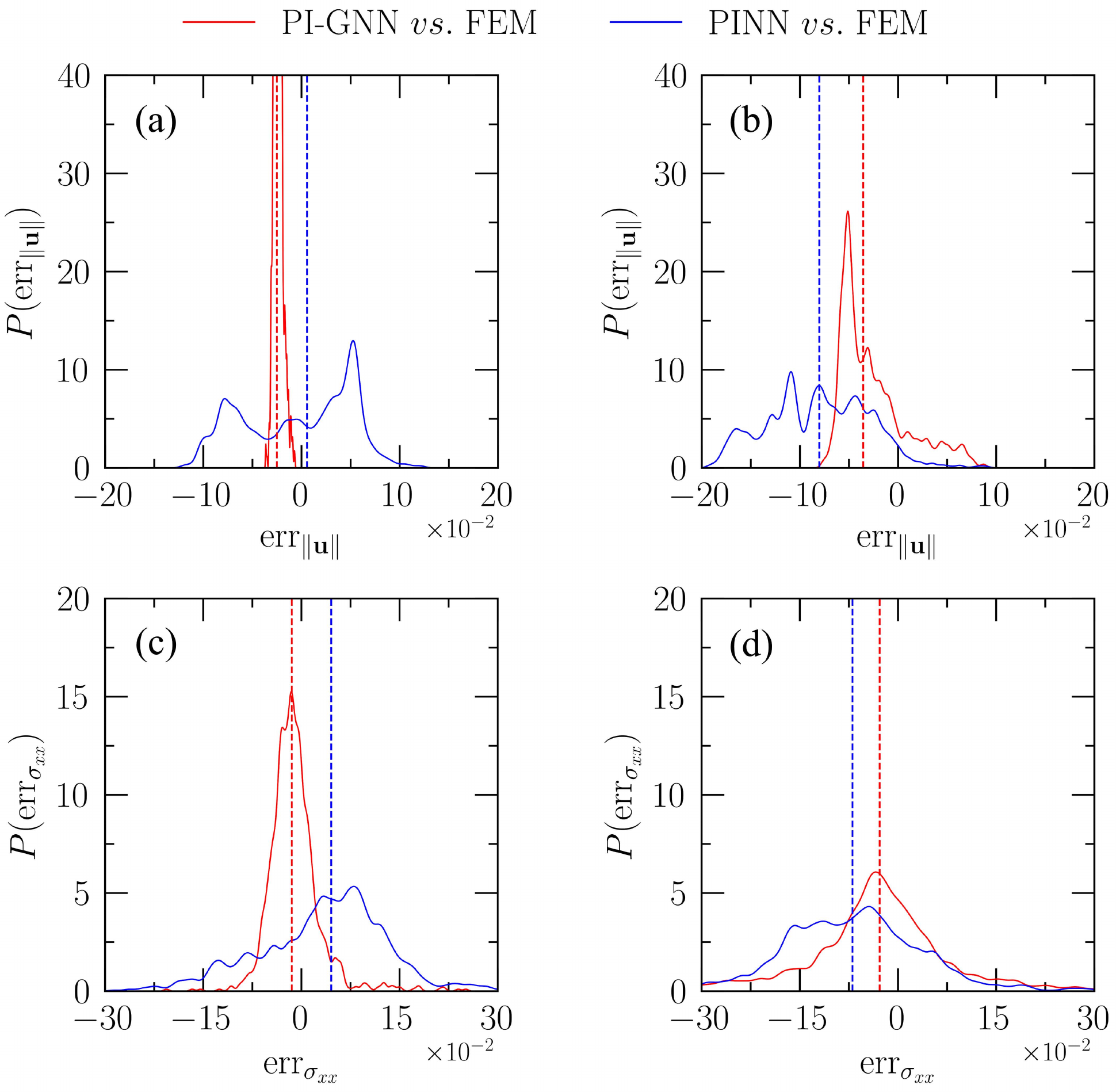}
    \caption{Spatial density of the element-wise error ($\text{err}_{\Vert\vb{u}\Vert}$ and $\text{err}_{\sigma_{xx}}$) for the PINN and \pignn predictions against the FE reference, where the dashed vertical lines are the median error for the respective solver. (a) \& (c) Single hole ($r = \SI{4}{\milli\meter}$). (b) \& (d) Asymmetric three-hole configuration. Nodes carrying less than $20\%$ of the peak stress are excluded to focus the comparison on the failure-relevant region.}

    \label{fig:disp_err_density}
\end{figure}

The \(\sigma_{xx}\) fields are compared in \Cref{fig:reference_multihole_configuration}. The close spacing between the holes concentrates the load within the narrow ligaments separating adjacent holes, with the peak \(\sigma_{xx}\) reaching \(\approx 5T\) at the upper edge of the largest hole. Both neural solvers reproduce the overall stress pattern qualitatively, but the quantitative metrics distinguish their accuracy (\Cref{tab:pinn_vs_pignn_scores}). The \pignn achieves \(R^2 \approx 0.96\) both globally and in the near-field, compared with \(\approx\ 0.94\) for the PINN, while its relative \(L^{2}_{\sigma_{xx}}\) error is approximately half that of PINN (\(7.57\%\) versus \(15.08\%\)). The \(\text{err}_{\sigma_{xx}}\) distribution (\Cref{fig:disp_err_density}d) show a similar trend to the single-hole example: the \pignn remains centered (median \(-0.68\%\), standard deviation of \(7.87\%\)), whereas the PINN exhibits a larger negative bias and broader distribution, with a median of \(-6.76\%\) and standard deviation of \(11.53\%\). Notably, the sign of the PINN bias changes between the two geometries, from positive in the single-hole example to negative in the three-hole example, whereas the \pignn remains close to unbiased in both examples. These results are consistent with the limitation reported by~\citet{Li2021}, who noted that although the energy-based model qualitatively captures the global stress distribution, its local prediction error remains appreciable in regions surrounding stress concentrations.

\begin{table}[H]
    \centering
    \small
    \setlength{\tabcolsep}{10pt} 
    \caption{Longitudinal stress (\(\sigma_{xx}\)) error metrics for the PINN and \pignn predictions against the FE reference on the unseen evaluation mesh, for the single-hole (\(r = \SI{4}{\milli\meter}\)) and three-hole configuration. The \(\text{err}_{\sigma_{xx}}\) distribution is summarized by its median and standard deviation across nodes that carry more than \(20\%\) of the peak stress.}
    \label{tab:pinn_vs_pignn_scores}
    
    \begin{tabular}{l *{2}{S[table-format=2.4] S[table-format=-1.4]}}
        \toprule \toprule
        & \multicolumn{2}{c}{Single hole} & \multicolumn{2}{c}{Three holes} \\
        \cmidrule(r){2-3} \cmidrule(l){4-5}
        Metric & {PINN} & {\pignn} & {PINN} & {\pignn} \\
        \midrule
        \(L^{2}_{\sigma_{xx}}\) (\%)  & 10.9425 &  5.0128 & 15.0781 &  7.5738 \\
        \(R^2\) (global)                        &  0.9686 &  0.9918 &  0.9384 &  0.9600 \\
        \(R^2\) (near-field)                    &  0.9715 &  0.9919 &  0.9419 &  0.9593 \\
        Std.\ dev.\ (\%)                        &  9.6903 &  5.2359 & 11.5258 &  7.8710 \\
        Median signed rel.\ error (\%)          &  5.0698 & -0.4116 & -6.7578 & -0.6793 \\
        \bottomrule \bottomrule
    \end{tabular}
\end{table}

Across the three configurations, the \pignn reproduces the analytical Kirsch solution for the single-hole example and achieves accuracy comparable to or higher than the PINN across the perforated-plate benchmarks, with the largest differences observed in the near-field regions surrounding the holes. The \pignn also maintains a median pointwise error close to zero across the considered configurations, whereas the PINN exhibits geometry-dependent bias. Both comparisons are performed under discretization transfer, with each model trained on its own sampling of the domain and evaluated on a separate unseen FE mesh using the same architecture and training settings. The \pignn is, therefore, evaluated by a forward pass on an unseen mesh graph. In contrast, the coordinate-based PINN is evaluated by querying its learned continuous map at new coordinates. The sharper near-field accuracy observed for the \pignn is consistent with its mesh-based energy formulation, in which the internal energy is evaluated element-wise using exact element areas rather than estimated through Monte Carlo collocation~\citep{Li2021}. 

\subsubsection{Two-phase composite}\label{sec:pinn_vs_pignn_inclusion}
\noindent The perforated plate benchmarks isolate geometric stress concentrations in a homogeneous material. Heterogeneous material systems introduce an additional challenge because the material properties are discontinuous across the interface. For a perfectly bonded bimaterial interface, the displacement field remains \(C^0\)-continuous and the traction is continuous, whereas the in-plane stress components exhibit finite discontinuities across the interface. An energy-based continuum micromechanics PINN is evaluated on the same two-phase benchmark to compare its accuracy with the proposed graph-based solver and to establish a reference error level for heterogeneous material systems \citep{henkes2022physics}. 

The benchmark consists of a square plate \(\Omega=[0,\aashay{l}]^2\) with \(l=\SI{2}{\milli\meter}\), containing a circular inclusion, \(\Omega_{\mathrm{inc}}=\mathcal{B}_a(\vb{x}_c)\), of radius \(a=\SI{0.4}{\milli\meter}\) centered at \(\vb{x}_c=(1,1)\,\si{\milli\meter}\), with the surrounding matrix occupying \(\Omega_{\mathrm{mat}}=\Omega\setminus\overline{\Omega}_{\mathrm{inc}}\). Both phases have \(\nu_{\mathrm{mat}}=\nu_{\mathrm{inc}}=0.4\), while \(E_{\mathrm{mat}}=\SI{1500}{\mega\pascal}\) and \(E_{\mathrm{inc}}=\SI{10000}{\mega\pascal}\), yielding a stiffness contrast of 6.67. A uniform tensile traction of \(T=\SI{0.025}{\mega\pascal}\) is applied to the right edge, with the remaining boundaries subjected to the roller-and-pin constraints. The geometry, material properties, and loading are identical to those of the reference study.

Unlike the deep-energy baseline, which shares the variational formulation of the \pignn, the reference model is a strong-form, mixed-variable PINN that predicts \((u_x, u_y, \sigma_{xx},\) \(\sigma_{yy}, \sigma_{xy})\). Equilibrium and constitutive consistency are enforced through collocated residuals rather than a single energy functional. The implementation follows the code released by~\citet{henkes2022physics}, retaining the published architecture and training settings: four hidden layers of 64 neurons with \(\tanh\) activations and \texttt{LeCun-uniform} initialization, a uniform \(128\times128\) collocation grid, and \(5000\) \texttt{BFGS} iterations, yielding a final training loss of \(3.41\times10^{-6}\). Material heterogeneity is represented through the hyperbolic-tangent regularization of the Lam\'e parameters across a transition band of width \(\delta=0.01\). 
\aashay{With further model related details in \Cref{app:sweep}.}

The \pignn uses the same architecture and hyperparameters from the previous verification example and is trained for \(50{,}000\) epochs using \texttt{Adam} with \(\eta=10^{-3}\). The two models, therefore, use different optimization formulations: a strong-form residual minimization with \texttt{BFGS} for the reference PINN and discrete energy minimization with \texttt{Adam} for the \pignn. The only modification to the released PINN implementation concerns the post-processing. Rather than evaluating the trained network on the \(128\times128\) collocation grid used during training, the network is evaluated at an independent set of points, allowing both neural models and the FE reference to be compared on a common discretization (3034 nodes and 5944 elements), which is not used during training. The mesh conforms to the material interface, \aashay{with 252 nodes on the interface.} This \aashay{interface-conforming discretization} permits the stress field to take distinct values on either side of the interface while retaining the required displacement continuity across the perfectly bonded interface. The \pignn, therefore, represents the material discontinuity explicitly through the conforming mesh topology, without introducing a smoothing or regularization layer at the interface. 

\begin{figure}[H]
    \centering
    \includegraphics[width=0.9\linewidth]{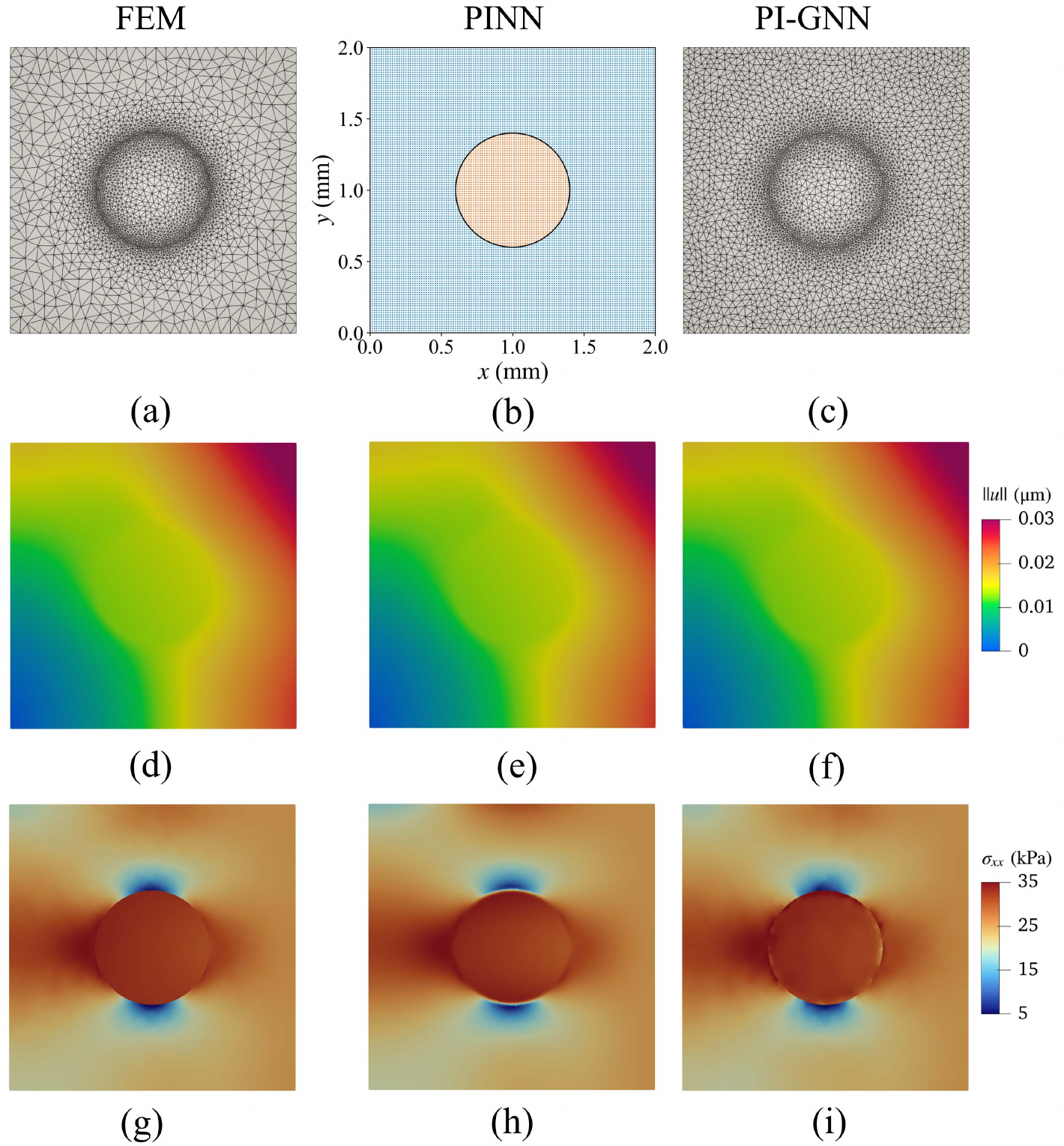}
    \caption{Discretizations and field predictions for the circular inclusion. 
    \textbf{Top row:} (a) Shared FE evaluation mesh ($3034$ nodes, $5944$ elements). (b) The uniform point cloud of \(128 \times 128\) grid on which the PINN of \citet{henkes2022physics} is trained. (c) \pignn training mesh ($3762$ nodes, $7316$ elements). The neural solvers are trained on (b) and (c), then evaluated on (a). 
    \textbf{Middle row:} Displacement magnitude ($\Vert\vb{u}\Vert$) fields evaluated on the FE reference mesh. (d) FE reference solution, (e) energy-based PINN prediction~\citep{henkes2022physics}, and (f) \pignn prediction. 
    \textbf{Bottom row:} Longitudinal stress ($\sigma_{xx}$) fields evaluated on the FE reference mesh. (g) FE reference solution, (h) \aashay{PINN prediction}, and (i) \pignn prediction.}
    \label{fig:pinn_pignn_inclusion_meshes}
\end{figure}

\Cref{fig:pinn_pignn_inclusion_meshes} compares the \(\Vert\vb{u}\Vert\) and \(\sigma_{xx}\) fields. Both models reproduce the characteristic load-transfer mechanism of a stiff inclusion embedded in a compliant matrix: the inclusion carries a large fraction of the applied load, resulting in elevated \(\sigma_{xx}\) within the inclusion and stress shielding in the surrounding matrix. The quantitative comparison shows comparable accuracy for the two models (see~\Cref{tab:pinn_vs_pignn_inclusion_l2,tab:pinn_vs_pignn_inclusion_r2}). \aashay{The relative \(L^2\) error in the displacement field is below \(1\%\) for both methods, with \(0.86\%\) for the PINN and \(0.17\%\) for the \pignn. For \(\sigma_{xx}\), the corresponding errors are \(3.02\%\) and \(3.78\%\), respectively, while the von Mises stress errors are \(2.57\%\) and \(1.85\%\). All predicted fields achieve \(R^2 \ge 0.95\) globally and within the near-field region. The transverse normal and shear stresses exhibit larger relative \(L^2\) errors, ranging from \(15\%\) to \(28\%\) for both models. These components carry little of the load under uniaxial tension; consequently, their larger relative errors primarily result from normalization by a small reference-field norm, while the corresponding absolute errors remain small (see~\Cref{sec:verification_summary}).} Overall, neither model exhibits a consistent accuracy advantage. This benchmark demonstrates that the proposed \pignn, without architecture or hyperparameter tuning, achieves accuracy comparable to that of a PINN specifically developed for heterogeneous continuum micromechanics.

Two differences nevertheless emerge in the distribution of the errors. The first concerns the near-field region. Restricting the evaluation to the 2668 nodes located within twice the inclusion radius of the interface, the \pignn exhibits improved \(R^2\) values for the displacement components, \(\sigma_{xx}\), and \(\sigma_\text{vM}\), whereas the PINN shows a slight reduction in \(R^2\) for all fields except \(\sigma_{xy}\). The signed errors show a similar trend: the \pignn remains centered near \(0.0\%\) for \(u_x\), \(\sigma_{xx}\), and \(\sigma_\text{vM}\), whereas the PINN exhibits biases of \(-0.86\%\) for \(u_x\). Although these differences are small, they indicate a modest advantage for the \pignn in resolving the fields near the material interface.

\begin{figure}[H]
    \centering
    \includegraphics[width=0.9\linewidth]{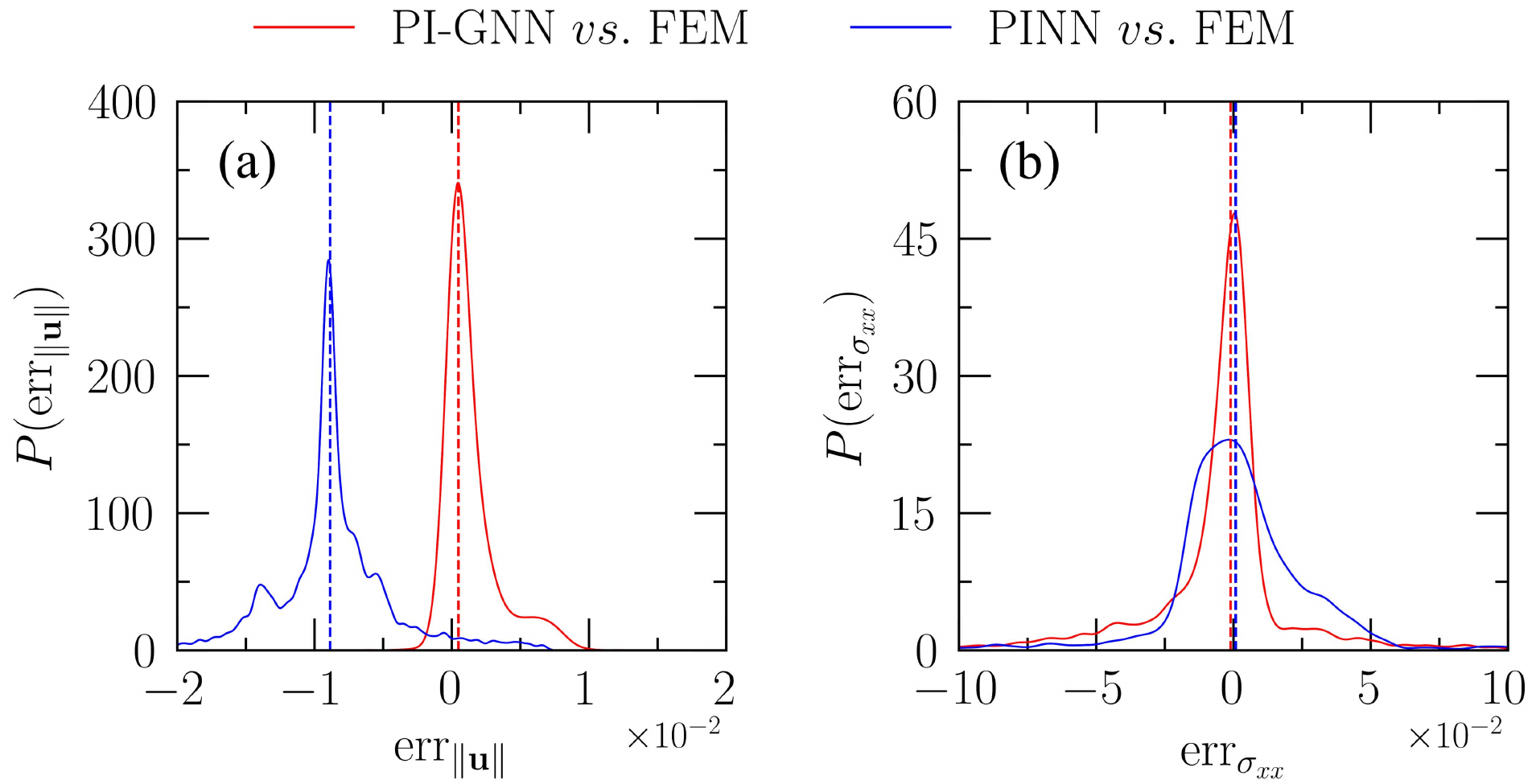}
    \caption{Spatial density of the element-wise relative error for the PINN and \pignn predictions against their respective FE references for the circular inclusion. (a) \(\text{err}_{||\vb{u}||}\). (b) \(\text{err}_{\sigma_{xx}}\). Nodes carrying less than \(20\%\) of the peak reference value are excluded, following~\Cref{eq:4.1}.}
    \label{fig:inclusion_err_density}
\end{figure}

The second difference concerns the interface's representation. Owing to the conforming mesh, the \pignn captures the stress discontinuity across \(\Gamma_{\mathrm{int}}\) directly through material assignment of adjacent elements. In contrast, the coordinate-based PINN represents the material interface through a hyperbolic-tangent transition band of width \(\delta=0.01\), producing a more diffuse transition in the predicted stress contours (see~\Cref{fig:pinn_pignn_inclusion_meshes}). The corresponding pointwise error distributions are shown in~\Cref{fig:inclusion_err_density}.

\begin{table}[H]
    \centering
    \small
    \setlength{\tabcolsep}{14pt}     
    \caption{\(L^2\) errors for the PINN and the \pignn, each against the FE solution, on the common evaluation mesh of the circular inclusion.} 
    \label{tab:pinn_vs_pignn_inclusion_l2}
    
    \begin{tabular}{l S[table-format=2.4] S[table-format=2.4]}
        \toprule \toprule 
        & \multicolumn{2}{c}{Relative Error \(L^{2}_\alpha\) (\%)} \\
        \cmidrule(l){2-3}
        Field & {PINN} & {\pignn} \\
        \midrule
        \(||\vb{u}||\) &  0.8622 &  0.1699 \\
        \(\sigma_{xx}\)    &  3.0223 &  3.7801 \\
        \(\sigma_{yy}\)    & 19.7799 & 27.4440 \\
        \(\sigma_{xy}\)    & 16.9558 & 14.5477 \\
        \(\sigma_\text{vM}\)    &  2.5672 &  1.8468 \\
        \bottomrule \bottomrule
    \end{tabular}
\end{table}

\begin{table}[H]
    \centering
    \small
    \setlength{\tabcolsep}{5pt}     
    \caption{Field-wise \(R^2\) scores and median signed relative errors. The near-field comprises the \(2668\) nodes within twice the inclusion radius of the interface.}
    \label{tab:pinn_vs_pignn_inclusion_r2}
    
    \begin{tabular}{l *{2}{S[table-format=1.4]} S[table-format=-1.4] *{2}{S[table-format=1.4]} S[table-format=-1.4]}
        \toprule \toprule
        & \multicolumn{3}{c}{PINN} & \multicolumn{3}{c}{\pignn} \\
        \cmidrule(r){2-4} \cmidrule(l){5-7}
        & \multicolumn{2}{c}{\(R^2\)} & {Median Error} & \multicolumn{2}{c}{\(R^2\)} & {Median Error} \\
        \cmidrule(r){2-3} \cmidrule(l){5-6}
        Field & {Global} & {Near-field} & {(\%)} & {Global} & {Near-field} & {(\%)} \\
        \midrule
        \(u_x\)          & 0.9989 & 0.9962 & -0.8570 & 0.9978 & 0.9999 &  0.1521 \\
        \(u_y\)          & 0.9988 & 0.9960 &  0.9318 & 0.9947 & 0.9999 &  -0.0951\\
        \(\sigma_{xx}\)  & 0.9820 & 0.9811 &  0.0921 & 0.9750 & 0.9855 &  0.1459 \\
        \(\sigma_{yy}\)  & 0.9740 & 0.9729 & -0.4006 & 0.9582 & 0.9571 &  1.5097 \\
        \(\sigma_{xy}\)  & 0.9731 & 0.9743 & -0.4342 & 0.9884 & 0.9882 &  0.0066 \\
        \(\sigma_\text{vM}\)  & 0.9796 & 0.9773 &  0.8186 & 0.9838 & 0.9929 & -0.0294 \\
        \bottomrule \bottomrule
    \end{tabular}
\end{table}

\noindent \textbf{Accuracy under increasing stiffness contrast:}
\amiya{The PINN and \pignn comparison with the FE reference solution described above considers a single stiffness contrast, \(E_{\mathrm{inc}}/E_{\mathrm{mat}}=6.67\). To assess the robustness of both solvers to increasing material contrast, we repeat the comparison across twenty-five logarithmically spaced contrasts spanning \(\Econ \in [10^{-2}, 10^{2}]\), placed reciprocally symmetrically about unit contrast, varying only \(E_{\mathrm{inc}}\). The matrix modulus, both Poisson's ratios, geometry, applied traction, evaluation mesh, and all architectural and optimization settings are kept fixed. In particular, the PINN retains its transition width \(\delta=0.01\), \(128\times128\) collocation grid, and \texttt{BFGS} iteration budget, with no problem-specific retuning. For each stiffness contrast, the solvers are trained from three independent weight initializations and evaluated on the same mesh using the area-weighted relative \(L^2\) (\Cref{eq:4.2}).}
\begin{figure}[H]
    \centering
    \subfloat[]{\includegraphics[width=0.47\linewidth]{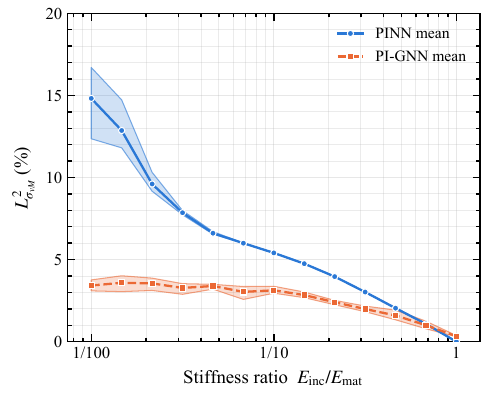}}\hspace{1em}
    \subfloat[]{\includegraphics[width=0.47\linewidth]{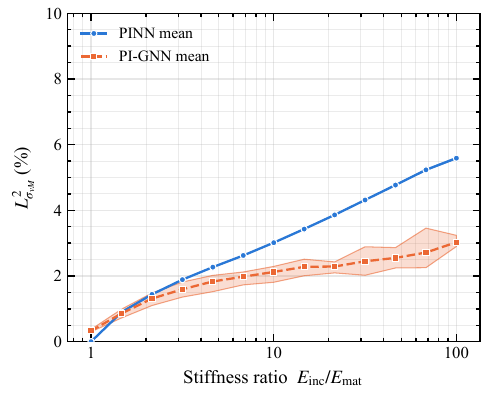}}
    \caption{\(L^{2}_{\sigma_\text{vM}}\) error against the FE reference as a function of the stiffness contrast, for the strong-form PINN~\citep{henkes2022physics} and the \pignn. Twenty-five logarithmically spaced contrasts spanning \(\Econ\ \in [10^{-2},10^2]\), placed reciprocally symmetrically about unit contrast. \aashay{Solid and dashed lines denote PINN and \pignn, respectively}. The lines show the mean across three random seeds, and the shaded regions show the seed-level range. The solvers are evaluated on an interface-conforming mesh (3034 nodes, 5944 elements) that is not used in \pignn training, and the area-weighted metric is used. (a) is the soft-inclusion branch and (b) the stiff branch.}
    \label{fig:stiffness_sweep}
\end{figure}
\amiya{The two solvers exhibit comparable accuracy for stiffness contrasts below approximately \(2\), with \(\sigma_{\mathrm{vM}}\) errors near \(2.3\%\) (see~\Cref{fig:stiffness_sweep}). The errors begin to separate beyond a contrast of approximately \(3\): the PINN error increases linearly with \(\log(E_{\mathrm{inc}}/E_{\mathrm{mat}})\), reaching \(5.58\%\) at a contrast of \(100\), whereas the \pignn remains below \(3.58\%\). The same trend is reflected in \(R^2\), with the PINN decreasing from a maximum of \(0.983\) to \(0.878\), while the \pignn remains near \(0.976\). The displacement error shows an even stronger separation, reaching \(3.10\%\) for the stiff-inclusion example and \(7.66\%\) for the soft-inclusion example at a contrast of \(100\), compared with \(0.20\%\) and \(0.49\%\), respectively, for the \pignn.}

\amiya{The \pignn is largely insensitive to whether the inclusion is stiffer or softer, whereas the PINN degradation is associated with the fixed interface regularization width \(\delta=0.01\). A contrast-dependent choice of \(\delta\) could reduce this degradation, but would introduce an additional problem-specific parameter. The comparison, therefore, does not claim superiority over an optimally tuned strong-form PINN. It rather demonstrates that the \pignn maintains a single configuration across the full contrast range without requiring an interface parameter to be adjusted. At a contrast of \(100\), the PINN error is 1.56 times the \pignn error for \(\sigma_{\mathrm{vM}}\). Repeating the evaluation on the \pignn training mesh gives \(3.07\%\), compared with \(3.58\%\) on the held-out mesh, indicating that the observed trend is not sensitive to the choice of evaluation mesh.}

At the reference contrast of \(6.67\), the two solvers achieve comparable accuracy at similar computational cost, despite using different optimization strategies. The \pignn completes \(50{,}000\) \texttt{Adam} epochs in \(\SI{2497}{\second}\), compared with \(\SI{1200}{\second}\) for the PINN after \(5000\) \texttt{BFGS} iterations on the same hardware. Thus, first-order optimization of the discrete potential energy achieves accuracy comparable to that of quasi-Newton optimization of the strong-form residuals on this benchmark. The memory requirements of \texttt{BFGS} also increase with the number of trainable parameters, which can become an important consideration for larger 3D and finite-strain problems. More importantly, the \pignn achieves this accuracy without problem-specific tuning. Geometry is introduced exclusively through the conforming mesh, while material heterogeneity enters through element-wise material properties. In contrast, coordinate-based PINNs require an appropriate collocation strategy and, for sharp material interfaces, a suitable interface regularization width. 

\subsubsection{Isolating the effect of message passing}\label{sec:message_passing}
\noindent \amiya{The comparisons above demonstrate an advantage over the collocation-based energy PINN of~\citet{Li2021} and, at high stiffness contrast, over the strong-form PINN of~\citet{henkes2022physics}. However, these comparisons do not isolate whether the observed advantage arises from the element-wise energy discretization or from message passing on the graph. To separate these effects, we construct a mesh-based PINN by modifying the existing \pignn while keeping its architecture and parameter count unchanged. Specifically, the neighborhood \(\mathcal{N}(i)\) in~\Cref{eq:message} is restricted to the node itself, eliminating information exchange between neighboring nodes. The resulting model, therefore, minimizes the same FE Ritz functional on the conforming mesh. Still, the mesh topology enters only through the assembly of the energy and does not influence the nodal trial field through message passing. This provides a direct comparison between mesh-based energy minimization with and without graph-based neighborhood interactions.}

\aashay{Two sweeps are performed on the two-phase inclusion problem. The first varies the mesh element size over nine values \(h \in [0.01, 0.05]\) at a fixed stiffness contrast; the second varies the contrast over seven values \(\Econ \in [10^{-2}, 10^{2}]\) on a fixed mesh. Each configuration is trained from 3 random seeds (initializations). The models use a reduced setting for this study \(17{,}250\) parameters (32 hidden neurons and 4 layers) and \(5000\) \texttt{Adam} epochs at \(\eta = 10^{-3}\). The smaller network size and limited training budget assess the convergence speed of the mesh-based PINN and \pignn models. At this fixed budget, the \pignn consistently achieves the lower von Mises error in both the mesh sweep and the stiffness-ratio sweep.}

\begin{figure}[H]
    \centering
    \subfloat[]{\includegraphics[width=0.45\linewidth]{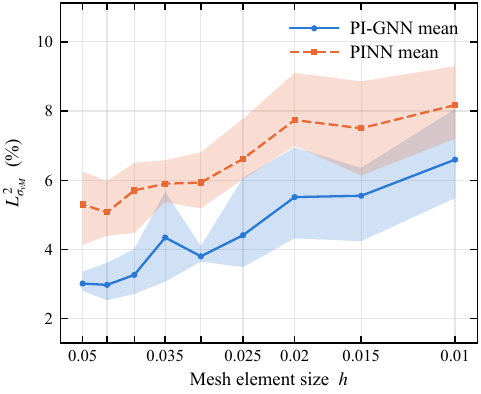}}\hspace{1em}
    \subfloat[]{\includegraphics[width=0.45\linewidth]{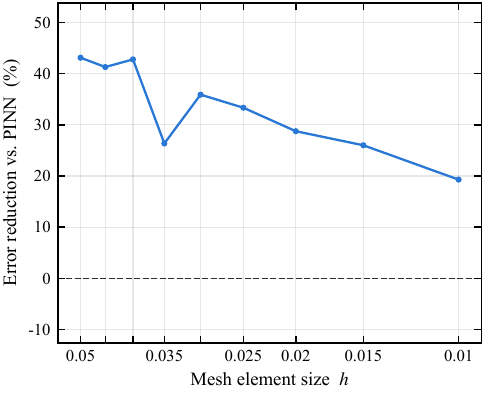}}\\
    \subfloat[]{\includegraphics[width=0.45\linewidth]{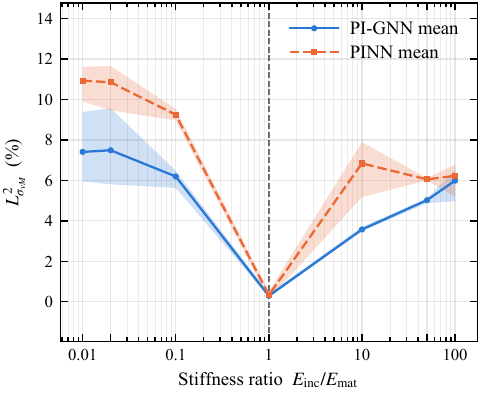}}\hspace{1em}
    \subfloat[]{\includegraphics[width=0.45\linewidth]{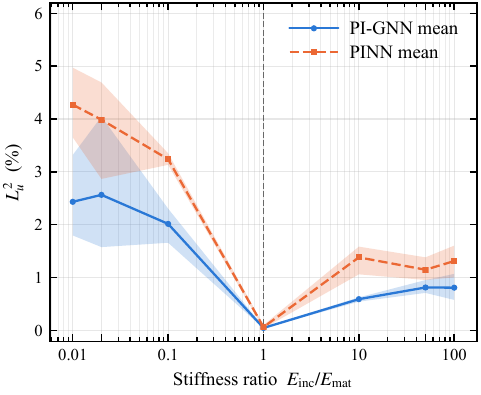}}
    \caption{Effect of message passing at a fixed parameter count. The mesh-based PINN is the \pignn with \(\mathcal{N}(i)\) restricted to the node itself. (a) \(L^{2}_{\sigma_\text{vM}}\) against the FE reference over the mesh sweep, and (b) the corresponding error reduction; \(h\) decreases to the right. (c) \(L^{2}_{\sigma_\text{vM}}\) over the stiffness-contrast sweep, and (d) \(L^{2}_{\text{u}}\) over the stiffness-contrast sweep ; the dashed vertical line marks the homogeneous limit \(\Econ = 1\). Lines show the mean across three random seeds, and shaded regions show the seed-level range.}
    \label{fig:message_passing}
\end{figure}

\aashay{Two trends indicate where the gain originates. The accuracy gain due to message passing increases as the mesh coarsens, from a 19\% error drop at \(h=0.01\) to \(43\%\) at \(h=0.05\) (see~\Cref{fig:message_passing}a \& b). With (\texttt{L = 4}) message-passing layers, a node aggregates over a graph neighborhood spanning four element widths, so a single forward pass covers a larger fraction of the inclusion when the mesh is coarse; as (h) decreases, the same four layers span a shorter physical distance and the margin narrows, although it remains at (19\%) on the finest mesh. The second trend is the stiffness contrast dependence (\Cref{fig:message_passing}c \& d). Here, the unit ratio refers to a homogeneous body, with no interface; both models yield errors near \(0.3\%\), and the gap closes. Away from the unit contrast, the gain in error due to message passing increases in displacement and stress fields. At this budget, the gain is therefore largest where the neighborhood exchange carries the most information, namely, across the material interface. Message passing incurs training time that is \(1.2\) to \(2.2\) times that of the ablated model at the same parameter count, because each epoch evaluates a message on every edge.}

\aashay{The models exhibit increasing error as \(h\) decreases, which reflects the fixed optimization budget rather than a loss of consistency. In the reduced setting (17,250 trainable parameters and 5000 \texttt{Adam} epochs), the optimization budget is held fixed while the number of nodal unknowns increases with mesh refinement. Consequently, the finer meshes require optimization over larger discrete systems and are less fully converged within the prescribed epoch budget. At the same time, mesh refinement yields a more accurate discrete approximation of the continuum solution, making the remaining optimization error more apparent when measured against the refined reference solution. Thus, \Cref{fig:message_passing}a reflects the optimization progress and convergence rate under a fixed computational budget, rather than the accuracy of fully converged models. A representative run with 30,000 epochs is provided in \ref{app:mp_converged}.
}

\subsection{Verification summary}\label{sec:verification_summary}
\noindent The verification suite demonstrates that a single \pignn architecture can accurately reproduce smooth elastic fields, localized stress concentrations, multi-feature stress interactions, and discontinuous constitutive responses at heterogeneous interfaces, without requiring labeled data and per-case tuning. \Cref{tab:verification_summary} summarizes the nodes and elements of the evaluation mesh, training configurations, and loss metrics for each example.  \amiya{The benchmark problems demonstrate that the proposed \pignn, without architecture or hyperparameter tuning, achieves accuracy comparable to that of a FEM/PINN.}

For the linear elastic problems, the internal strain energy \(U\) and external work \(W_{\text{ext}}\) satisfy Clapeyron's theorem (\(2U = W_{\text{ext}}\)) at equilibrium. Consequently, the total potential energy \(\Pi = U - W_{\text{ext}} = -U\), and the converged solution is expected to satisfy \(\Pi / W_{\text{ext}} = -0.5\). The convergence history is monitored to verify stabilization at this ratio. The converged \pignn energy is compared with the FE value using the relative difference \((\Pi_{\mathrm{PI\text{-}GNN}} - \Pi_{\mathrm{FEM}})/\lvert\Pi_{\mathrm{FEM}}\rvert\). 
\begin{table}[H]
\centering
\caption{Summary of evaluation mesh configuration, training, and metrics across the verification examples.}
\label{tab:verification_summary}
\begin{tabular}{lcccccccc}
\hline \hline
 & \multicolumn{2}{c}{Training mesh} & \multicolumn{2}{c}{FE evaluation mesh} & & & & \\
\cline{2-3} \cline{4-5}
Problem & Nodes & Elements & Nodes & Elements & Epochs & \(\eta\) & $\Pi/W_{\mathrm{ext}}$ & {$err_{\Pi}$} \\
\hline
Single-hole  & 4368 & 8221 & 1865 & 3420 & 30,000 & 0.001 & -0.49 & +0.40 \\
Three-hole  & 4455 & 8394 & 3158 & 5837 & 30,000 & 0.001 & -0.50 & +1.97 \\
Inclusion & 3762 & 7316 & 3034 & 5944 & 50,000 & 0.001 & -0.50 & +0.13 \\
\hline \hline
\end{tabular}
\end{table}

\noindent \textbf{Remark:} \aashay{In each example, the non-dominant stress components, such as \(\sigma_{yy}\) and \(\sigma_{xy}\) under uniaxial tension, carry a small fraction of the load and consequently exhibit larger relative \(L^2\) errors, even though the dominant stress component and \(\sigma_\mathrm{vM}\) are closely reproduced. This difference is primarily a consequence of the component-wise normalization in \Cref{eq:4.2}, where \(L^{2}_{\alpha}\) is divided by the norm of the corresponding reference component. Thus, for components of small magnitude, even a small absolute error can result in a comparatively large relative error. The nodal exclusion in \Cref{eq:4.1} does not alter this normalization, since it is applied independently to each stress component.} 

\aashay{The small magnitude of these errors is further supported by the close recovery of the von Mises stress, which differs from the reference by at most 1.85\%. An error in a non-dominant component comparable in magnitude to the dominant stress would generally produce a substantially larger deviation in the resulting von Mises field and is therefore inconsistent with the observed agreement. The coefficients of determination for these components also remain above 0.95, although \(R^2\) provides limited independent evidence in this setting (see \Cref{tab:pinn_vs_pignn_inclusion_r2}). For components with means close to zero, the normalization in \(R^2\) and the relative \(L^2\) error are governed by quantities of similar magnitude, making the two measures strongly related. For components that are identically zero by construction, such as the in-plane normal stresses under pure torsion, the reference variance vanishes, and \(R^2\) is consequently not informative. The same normalization argument applies to the non-dominant components in the remaining examples.}

\section{Numerical applications}\label{sec:applications}
\noindent This section applies the verified framework to three challenging problems for which no closed-form solution exists. In each example, a converged FE solution computed on the same conforming mesh serves as the reference, so the reported errors quantify only the approximation introduced by the \pignn. The applications extend the verification benchmarks in three complementary directions. The first considers re-entrant inclusions to assess the resolution of stress concentrations induced by concave corners. The second and third investigate a 3D reinforced cube and rod subjected to torsion, respectively. All problems are trained using the same network architecture and hyperparameters similar to \Cref{sec:verification}, with a constant learning rate of $\eta=1\times10^{-3}$ for 100,000 epochs. The training and evaluation are performed on the same conforming mesh.

\subsection{Material Discontinuities}\label{sec:irregular-inclusion}
\noindent A material discontinuity is introduced by embedding a stiff inclusion within a compliant matrix to assess the \pignn to resolve heterogeneous constitutive behavior across bonded material interfaces. The geometry and material properties are adopted from~\citet{Kamali2023}. Following the reference study, the surrounding matrix consists of a \(10\%\) w/v PVA hydrogel (\(E_{\mathrm{mat}}=\SI{1500}{\pascal}\), \(\nu_{\mathrm{mat}}=0.45\)), while the inclusion is a stiffer \(15\%\) w/v PVA hydrogel (\(E_{\mathrm{inc}}=\SI{5000}{\pascal}\), \(\nu_{\mathrm{inc}}=0.35\)), resulting in a stiffness contrast of \(E_{\mathrm{inc}}/E_{\mathrm{mat}}\approx3.33\). Both phases, matrix and inclusion, are modeled as Neo-Hookean solids. 

\subsubsection{Re-entrant inclusion}\label{sec:re-entrant_tension}
\noindent The circular inclusion example is extended to a non-convex inclusion geometry to evaluate the \pignn to resolve the steep stress gradients introduced by re-entrant corners. Unlike the smooth circular interface, the reversal of interface curvature at the re-entrant lobes produces localized stress concentrations that are challenging to capture. The geometry consists of a square hydrogel layer of side $l=\SI{2}{\milli\meter}$ containing a five-lobed re-entrant inclusion positioned slightly off the plate center, eliminating symmetry about both coordinate axes. The same Dirichlet boundary conditions are imposed, together with a uniform tensile traction of $T=\SI{50}{\pascal}$ applied to the right edge. The conforming FE mesh consists of 5088 nodes and 9966 linear triangular elements, with local refinement along the re-entrant interface. The union of the refined regions defines the near-field domain used for the localized error analysis.

\begin{figure}[H]
    \centering
    \includegraphics[width=0.9\linewidth]{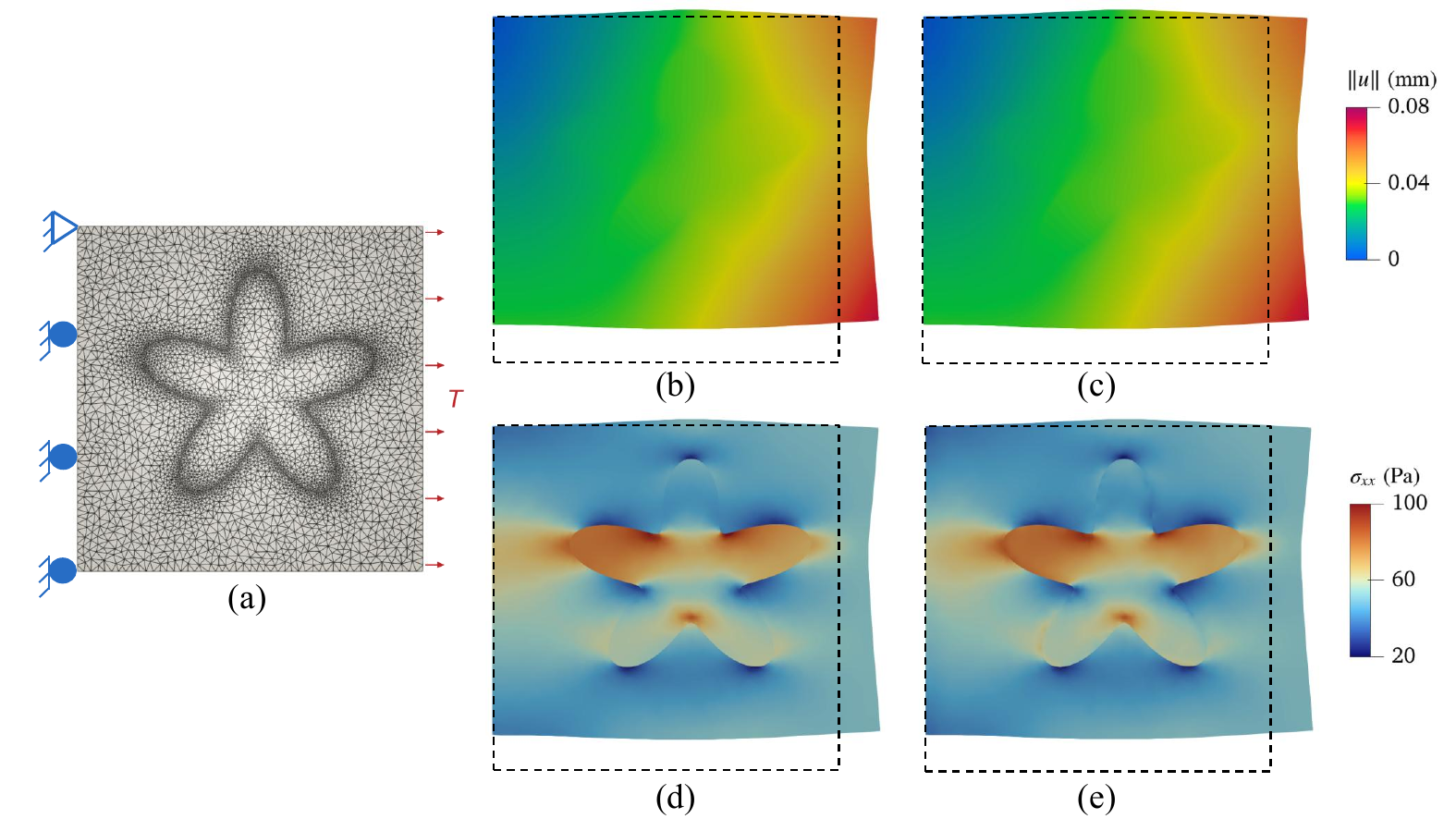}
    \caption{Stiff re-entrant inclusion in a soft hydrogel layer. (a) Reference configuration showing a conforming graph. (b) Deformed FEM configuration showing $||u||$; (c) Deformed \pignn prediction showing $||u||$; (d) Deformed FEM configuration showing $\sigma_{xx}$. (e) Deformed \pignn prediction showing $\sigma_{xx}$. Both stress panels use the same color scale in \SI{}{\pascal}. The deformed configurations are scaled by a factor of 5.}
    \label{fig:petal_3panel_sigma_xx}
\end{figure}

\Cref{fig:petal_3panel_sigma_xx} presents the stiffer inclusion carries the elevated stress, reaching $\approx \SI{100}{\pascal}$ along the lobes aligned with the loading direction. In contrast, the matrix immediately above and below the inclusion is shielded and relaxes toward $\SI{40}{\pascal}$. The \pignn reproduces the position and magnitude of every lobe, including the steep gradients at the concave notches between adjacent lobes. 

The displacement components achieve $R^2 \geq 0.99$ globally and in the near field, indicating that the irregular interface does not degrade the displacement solution (\Cref{tab:error_metrics_multi_inclusion}). The signed median errors for $u_x$ and $u_y$ are $-0.7161\%$ and $+0.6681\%$, indicating no appreciable bias in either component. The dominant stress $\sigma_{xx}$ and the von Mises stress achieve $R^2 \geq 0.98$, with near-field values marginally exceeding their global counterparts, indicating consistently accurate recovery of the primary load-carrying fields.
\begin{table}[H]
    \centering
    \small
    \setlength{\tabcolsep}{10pt}     
    \caption{Error metrics comparing the \pignn predictions with FEM reference for re-entrant inclusion.}
    \label{tab:error_metrics_multi_inclusion}
    
    \begin{tabular}{l S[table-format=1.4] S[table-format=-1.4] S[table-format=1.4] S[table-format=1.4]}
        \toprule \toprule
        & & & \multicolumn{2}{c}{\(R^2\)} \\
        \cmidrule(l){4-5}
        Field & {\(L^{2}_\alpha\) (\%)} & {Median \(\text{err}_{\alpha}\) (\%)} & {Global} & {Near-field} \\
        \midrule
        \(u_x\)         & 0.7108 & -0.7161 & 0.9956 & 0.9994 \\
        \(u_y\)         & 0.6211 &  0.6681 & 0.9971 & 0.9993 \\
        \(\sigma_{xx}\) & 2.7457 & -0.4408 & 0.9831 & 0.9908 \\
        \(\sigma_\text{vM}\) & 2.2346 & -0.2898 & 0.9852 & 0.9935 \\
        \bottomrule \bottomrule
    \end{tabular}
\end{table}

\begin{figure}[H]
    \centering
    \includegraphics[width=0.9\linewidth]{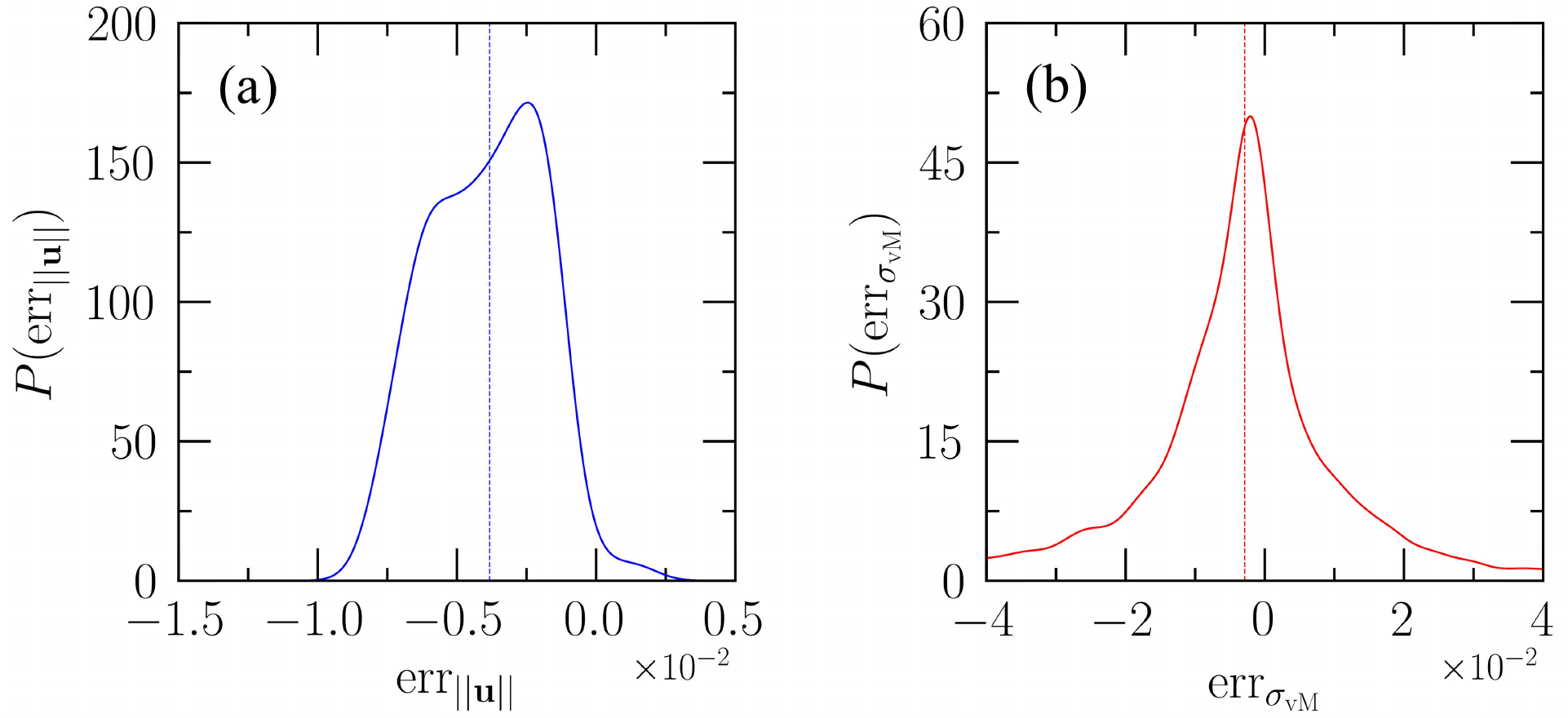}
    \caption{Relative error densities comparing \pignn and the FEM reference for the re-entrant inclusion. (a) Displacement magnitude. (b) von Mises Stress.}
    \label{fig:irregular_inc_errors}
\end{figure}

The error-density plots show how the residuals are distributed across the domain (\Cref{fig:irregular_inc_errors}). The displacement magnitude error is unimodal and tightly concentrated, with its median near $\text{err}_{||\vb{u}||} \approx -0.3\%$ and most of the probability mass between $-1.0\%$ and $+0.3\%$ (see~\Cref{fig:irregular_inc_errors}a). The bias is therefore negligible and marginally negative. The von Mises error is centered near $\text{err}_{\sigma_\text{vM}} \approx -0.3\%$, with a roughly symmetric core inside $\pm 3\%$ (see~\Cref{fig:irregular_inc_errors}b). Further, the displacement field is recovered with $R^2 > 0.99$, and the stress concentrations at the concave notches are resolved without a local penalty, since the near-field $R^2$ exceeds the global value for every field. The error metrics indicate that the \pignn framework handles circular and re-entrant inclusions without bias, and that accuracy does not deteriorate at the interface, where the gradients are steepest. 

\subsubsection{Re-entrant inclusion with combined normal-shear loading}\label{sec:ang_plate}
\begin{figure}[ht]
    \centering
    \includegraphics[width=0.9\linewidth]{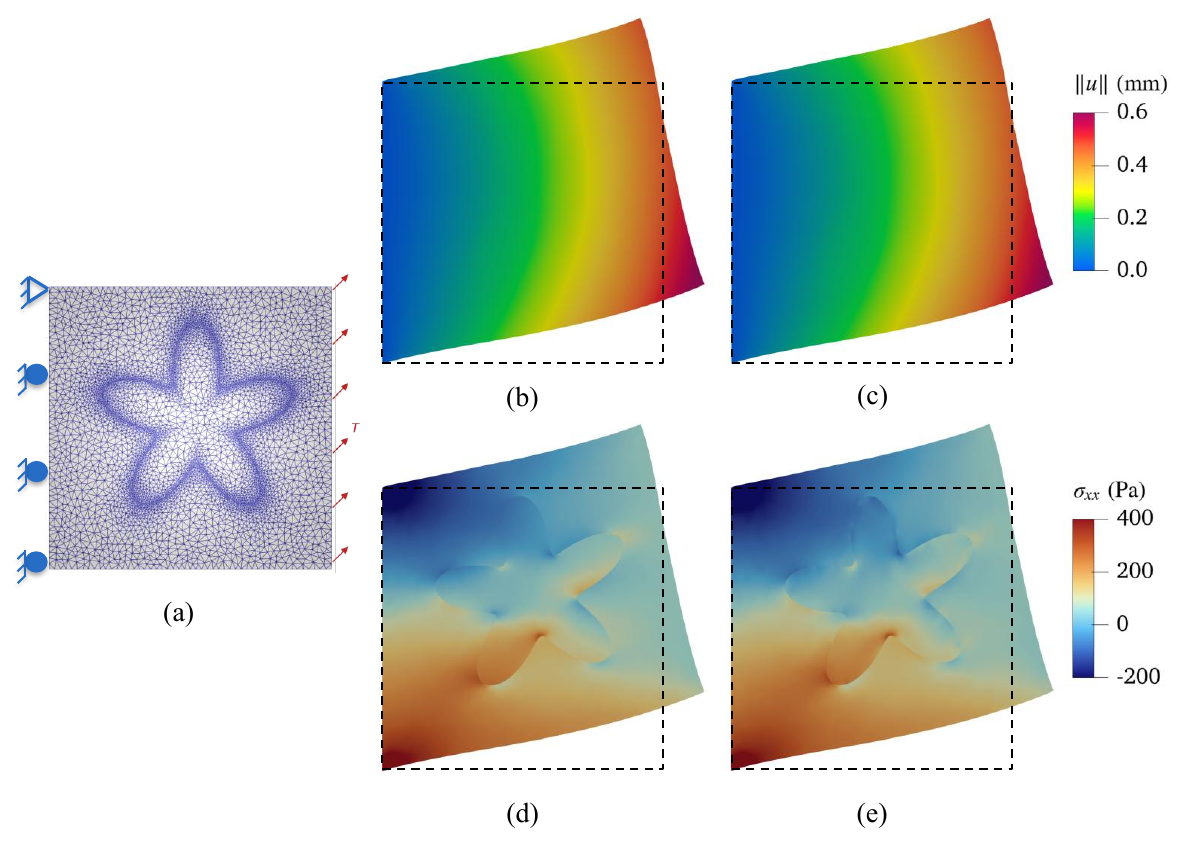}
    \caption{Stiff re-entrant inclusion in a soft hydrogel layer under normal-shear loading. The stress fields are shown using the same scaled colormap. (a) Reference configuration showing conforming graph and updated boundary conditions. (b) Deformed FEM configuration showing $||u||$; (c) Deformed \pignn prediction showing $||u||$; (d) Deformed FEM configuration showing $\sigma_{xx}$. (e) Deformed \pignn prediction showing $\sigma_{xx}$. Both stress panels use the same color scale in \SI{}{\pascal}. The deformed configurations are scaled by a factor of 5.}
    \label{fig:petal_angular_sigma_xx}
\end{figure}

\noindent The re-entrant inclusion problem is subjected to combined normal and shear loading, unlike the uniaxial example, where loading is carried entirely by the longitudinal normal stress, the applied traction $T=\SI{100}{\pascal}$ is now inclined at \SI{45}{\degree} to the \(x\)-axis. The resulting loading produces equal normal and transverse traction components, resulting in a coupled normal-shear stress state. To ensure static equilibrium, the left boundary is fully clamped ($u_x=u_y=0$) instead of being supported by a roller with a single corner pin. The transverse component of the applied traction introduces a net vertical force and bending moment that cannot be balanced by the roller constraint, making a fully fixed boundary essential for a well-posed problem (see~\Cref{fig:petal_angular_sigma_xx}a). 

\Cref{fig:petal_angular_sigma_xx}b \& c compare the reference and \pignn predictions of the $\sigma_{xx}$ field on the deformed configuration, which exhibits both rigid-body rotation and finite stretch. In contrast to the pure tension, the inclined traction acting against the clamped boundary yields a bending response. The transverse component of the applied load generates a bending moment at the fixed edge, causing $\sigma_{xx}$ to vary from compression in the upper region (approximately $\SI{-200}{\pascal}$) to tension in the lower region (approximately $\SI{400}{\pascal}$), superposed on the axial tensile stress. As in the pure-tension example, the stiff re-entrant inclusion attracts load into the lobes aligned with the resultant traction, while the concave notches between adjacent lobes continue to exhibit the steepest stress gradients. 

The quantitative error metrics are summarized in~\Cref{tab:error_metrics_ang_inclusion}. The displacement field is recovered with good accuracy, with $u_x$ and $u_y$ achieving $R^2>0.99$ globally and in the near-field region. Thus, neither the inclined loading nor the fully clamped boundary adversely affects the displacement prediction. The stress metrics provide a more demanding assessment. Under the uniaxial loading, the shear stress $\sigma_{xy}$ was a low-magnitude component. Rotating the applied traction to \SI{45}{\degree} promotes $\sigma_{xy}$ to a primary load-carrying stress, yielding a relative error of $L^{2}_\alpha=6.3\%$, while maintaining $R^2>0.97$ for all stress components. The relative error of the dominant stress $\sigma_{xx}$ increases from $2.8\%$ to $6.8\%$, reflecting the transition from a predominantly tensile field to one containing both tensile and compressive regions. Nevertheless, its coefficient of determination remains high ($R^2 \ge 0.99$). 

\begin{figure}[H]
    \centering
    \includegraphics[width = 0.9\linewidth]{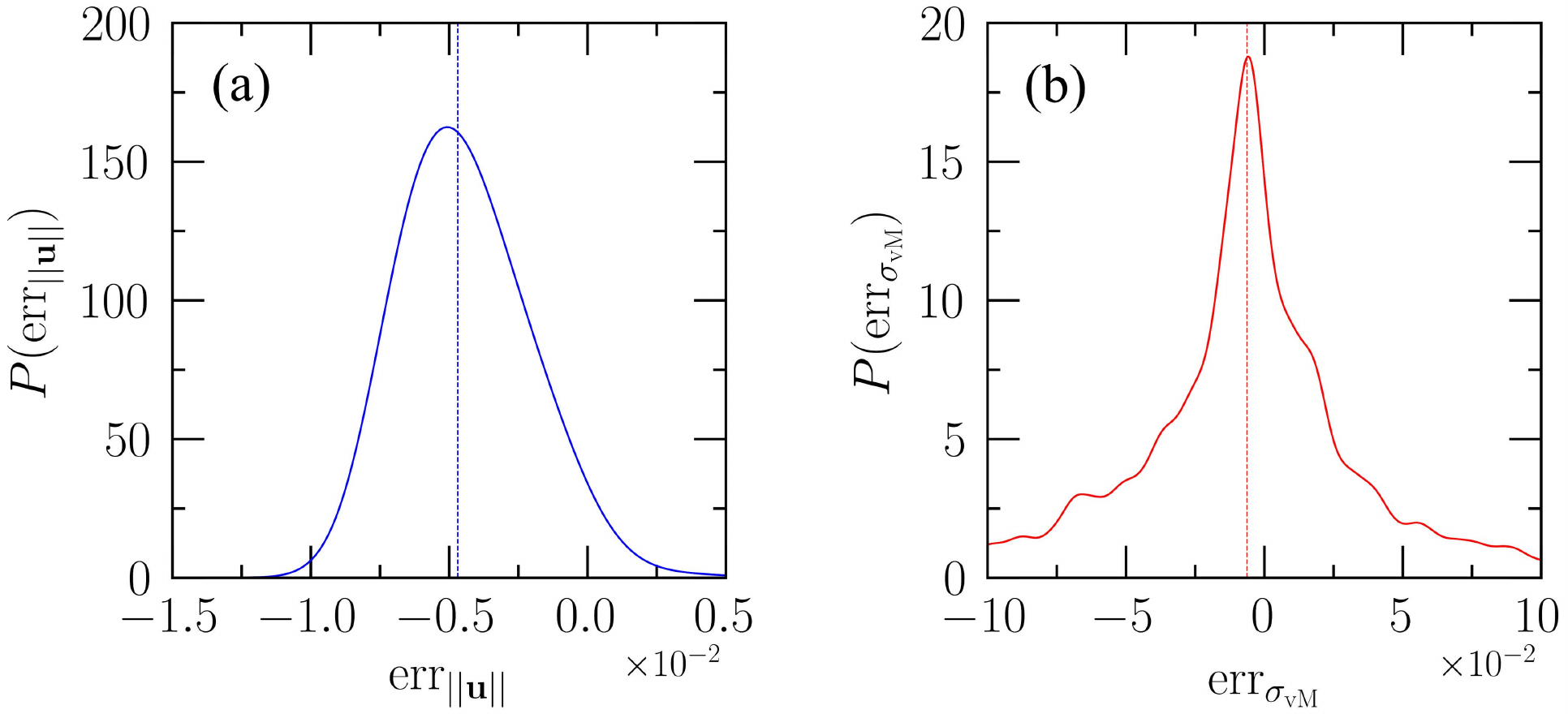}
    \caption{Relative error densities comparing \pignn and the FEM reference for the normal-shear loaded re-entrant inclusion. (a) Displacement magnitude. (b) von Mises Stress.}
    \label{fig:ang_inc_errors}
\end{figure}

\aashay{Compared to the uniaxial example, the $\sigma{vM}$ exhibits a moderate increase in relative error, from \(2.2\%\) to \(4.9\%\), consistent with its positive-definite nature.} Overall, the stress fields are recovered with high fidelity despite the increased complexity of the multiaxial stress state. \Cref{fig:ang_inc_errors} shows the pointwise error distributions, providing further insights. The displacement-magnitude error exhibits a narrow distribution with a median of approximately \(-0.5\%\), indicating a small systematic underprediction of the displacement magnitude (see \Cref{fig:ang_inc_errors}a). In contrast, the relative error of the $\sigma{vM}$ is centered close to zero (median $\approx -0.8\%$), with most values confined within $\pm6\%$ (see \Cref{fig:ang_inc_errors}b). These distributions indicate that the displacement and stress predictions remain essentially unbiased under combined loading.
\begin{table}[H]
    \centering
    \small
    \setlength{\tabcolsep}{10pt}    
    \caption{Error metrics comparing the \pignn predictions and FEM reference for re-entrant inclusion with normal-shear loading.}
    \label{tab:error_metrics_ang_inclusion}
    \begin{tabular}{l S[table-format=1.4] S[table-format=-1.4] S[table-format=1.4] S[table-format=1.4]}
        \toprule \toprule
        & & & \multicolumn{2}{c}{\(R^2\)} \\
        \cmidrule(l){4-5}
        Field & {\(L^{2}_\alpha\) (\%)} & {Median \(\text{err}_{\alpha}\) (\%)} & {Global} & {Near-field} \\
        \midrule
        \(u_x\)         & 1.3539 & -2.7129 & 0.9967 & 0.9997 \\
        \(u_y\)         & 0.5690 & -0.4687 & 0.9940 & 0.9998 \\
        \(\sigma_{xx}\) & 6.8704 & -5.9344 & 0.9907 & 0.9905 \\
        \(\sigma_{xy}\) & 6.2773 & -1.6033 & 0.9776 & 0.9736 \\
        \(\sigma_\text{vM}\) & 4.8927 & -0.8019 & 0.9776 & 0.9784 \\
        \bottomrule \bottomrule
    \end{tabular}
\end{table}

\subsection{Composite cube under torsion}\label{sec:cube_torsion}
\noindent We further consider a 3D composite cube under torsion to assess the generalizability of the proposed GNN and energy-based formulation beyond 2D problems. In this example, the geometry consists of a rubber cube (side \(\SI{20}{\milli\meter}\)) containing a concentric nylon core of square cross-section, \(\SI{8}{\milli\meter}\) on a side, extending over the full \(\SI{20}{\milli\meter}\) length along the \(z\)-axis. The rubber matrix has shear and bulk moduli \(\mu_{\mathrm{mat}}=\SI{7.41}{\mega\pascal}\) and \(\kappa_{\mathrm{mat}}=\SI{22.22}{\mega\pascal}\), corresponding to \(\nu_{\mathrm{mat}}\approx0.35\). The nylon core is assigned \(E_{\mathrm{inc}}=\SI{2000}{\mega\pascal}\) and \(\nu_{\mathrm{inc}}=0.3\), giving a stiffness ratio of \(E_{\mathrm{inc}}/E_{\mathrm{mat}}\approx100\). The face \(z=z_{\min}\) is fully fixed, \(\mathbf{u}=\mathbf{0}\), while the opposite face \(z=z_{\max}\) is subjected to the torsional traction \(\vb{t}=\tau_0(-y,x,0)\), producing a pure twisting moment about the cube axis. For \(\tau_0=\SI{0.15}{\mega\pascal}\), the resulting torque is \(\SI{4000}{\newton\milli\meter}\).

The conforming mesh contains \(9261\) nodes and \(48000\) tetrahedral elements. The \(d\)-dimensional architecture is instantiated with the feature vector
\(\vb{f}_i=[\hat{X}_i,\hat{Y}_i,\hat{Z}_i,m_i,\chi_i]^{\mathsf{T}}\in\mathbb{R}^5\)
and decoder output \(\tilde{\vb{u}}_i\in\mathbb{R}^3\). The Dirichlet mask is applied component-wise on the fixed face. The discrete energy is assembled element-wise over the tetrahedra, using the element volume \(V_e\) for the internal energy and integrating the torsional traction over the triangular facets of the loaded face, with facet area \(S_f\). Here, a uniform grid whose nodes align with the nylon core's faces is used, in contrast to the adaptive mesh used in previous examples. 

\begin{figure}[ht]
    \centering
    \includegraphics[width=0.9\linewidth]{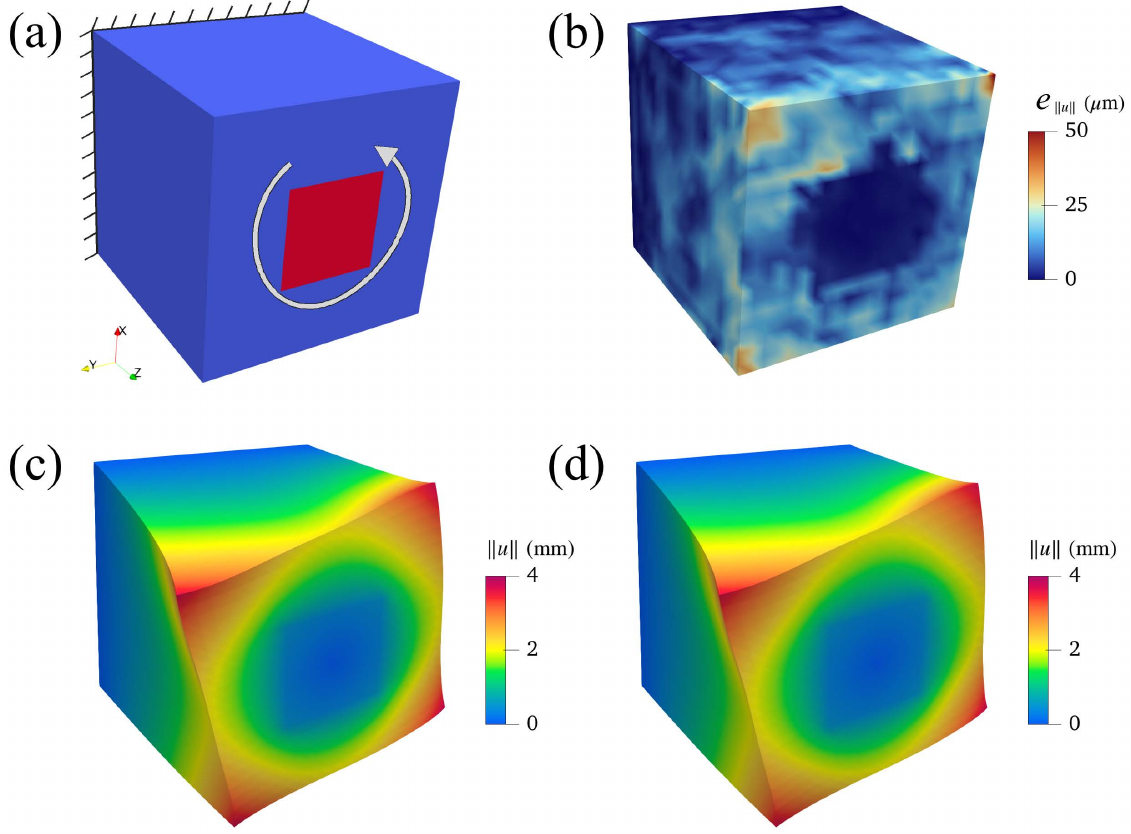}
    \caption{Hyperelastic composite cube under torsion. (a) Reference configuration with the applied torque. (b) Absolute displacement error \((e_{\lVert\vb{u}\rVert})\) between the reference FE solution and the \pignn prediction, mapped onto the cube faces. (c) Deformed configuration obtained from the FE displacement field \((\vb{u})\). (d) Deformed configuration from the \pignn prediction. Both deformed configurations are shown on the same 1:1 scale, with dimensions in \(\si{\milli\meter}\).}
    \label{fig:cube_disp}
\end{figure}

\Cref{fig:cube_disp} compares the deformed configurations for the \pignn and reference solver. The fixed face holds the cube in place while the loaded face rotates about the \(z\)-axis. The $\lVert\vb{u}\rVert$ grows from near zero at the center of each cross-section toward the corners, reaching $\approx \SI{4}{\milli\meter}$ at the free-face corners (see \Cref{fig:cube_disp}c \& d). The \pignn twist field (\Cref{fig:cube_disp}d) is visually indistinguishable from the FEM reference (\Cref{fig:cube_disp}c). The absolute error of displacement magnitude $e_{\lVert\vb{u}\rVert}$ mapped onto the cube faces stays below $\approx \SI{0.05}{\milli\meter}$ over most of the volume with $\approx \SI{0.1}{\milli\meter}$ near the corners (see~\Cref{fig:cube_disp}b). The stiffness contrast between the nylon core and the surrounding rubber results in a distinct local error pattern.

\Cref{tab:cube_error_stats} summarizes the error statistics for the torsion problem, indicating that the solution is resolved with accuracy comparable to the 2D inclusion examples. The in-plane displacements $u_x$ and $u_y$, which dominate the twist, and the axial component $u_z$, tied to warping, reach $R^2 \geq 0.99$. The stress field is shear-dominated, as expected under torsion: the two load-carrying shear components $\sigma_{xz}$, $\sigma_{yz}$, and invariant $\sigma_{vM}$ reach $R^2 \geq 0.99$, all with signed median errors below \(0.5\%\).
\begin{table}[ht]
    \centering
    \small
    \setlength{\tabcolsep}{14pt}     
    \caption{Summary of error statistics for the composite cube.}
    \label{tab:cube_error_stats}
    
    \begin{tabular}{l S[table-format=1.4] S[table-format=1.4] S[table-format=1.4]}
        \toprule \toprule
        Field & {\(L^{2}_\alpha\) (\%)} & {Median \(\text{err}_{\alpha}\) (\%)} & {\(R^2\) (global)} \\
        \midrule
        \(u_x\)         & 0.7651 & 0.2147 & 0.9999 \\
        \(u_y\)         & 0.8413 & 0.2990 & 0.9999 \\
        \(\sigma_{xz}\) & 4.2821 & 0.2281 & 0.9993 \\
        \(\sigma_{yz}\) & 4.2545 & 0.2269 & 0.9993 \\
        \(\sigma_\text{vM}\) & 3.1946 & 0.4636 & 0.9995 \\
        \bottomrule \bottomrule
    \end{tabular}
\end{table}

\begin{figure}[ht]
    \centering
    \includegraphics[width = 0.9\linewidth]{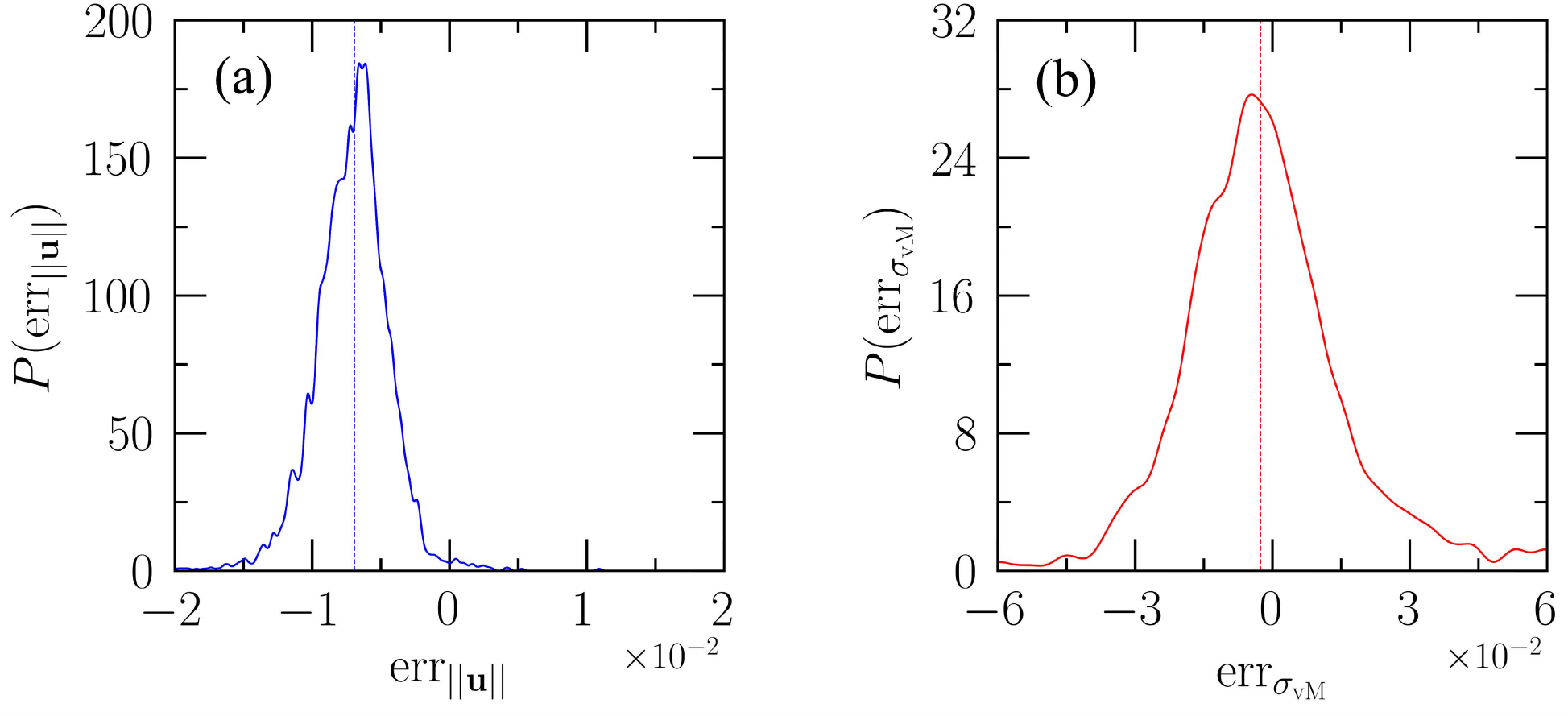}
    \caption{Comparing 3D \pignn against the FEM reference using \(err_{f}\) density. (a) Displacement magnitude. (b) von Mises Stress.}
    \label{fig:cube_errors}
\end{figure}

The error-density plots are consistent with these statistics (see~\Cref{fig:cube_errors}). The $\text{err}_{\lVert\vb{u}\rVert}$ (a) is unimodal with its median near \(-2.5\%\), so the network slightly under-predicts the twist amplitude. The von Mises error is broader and roughly symmetric about a small negative median \((\approx -0.5\%)\), with most of the probability mass inside $\pm 4\%$ (see~\Cref{fig:cube_errors}b). The spread is wider than the $\pm 2\%$ core seen for the re-entrant inclusion, reflecting the noisier stress recovery from a denser 3D graph (refer~\cref{sec:irregular-inclusion}). The dominant shear stresses and the von Mises invariant are all recovered with $R^2 > 0.99$, and the twist field matches the reference to within a couple of percentage points.

\subsection{Composite rod with a re-entrant inclusion under torsion}\label{sec:rod_torsion}
\begin{figure}[ht]
    \centering
    \includegraphics[width=0.9\linewidth]{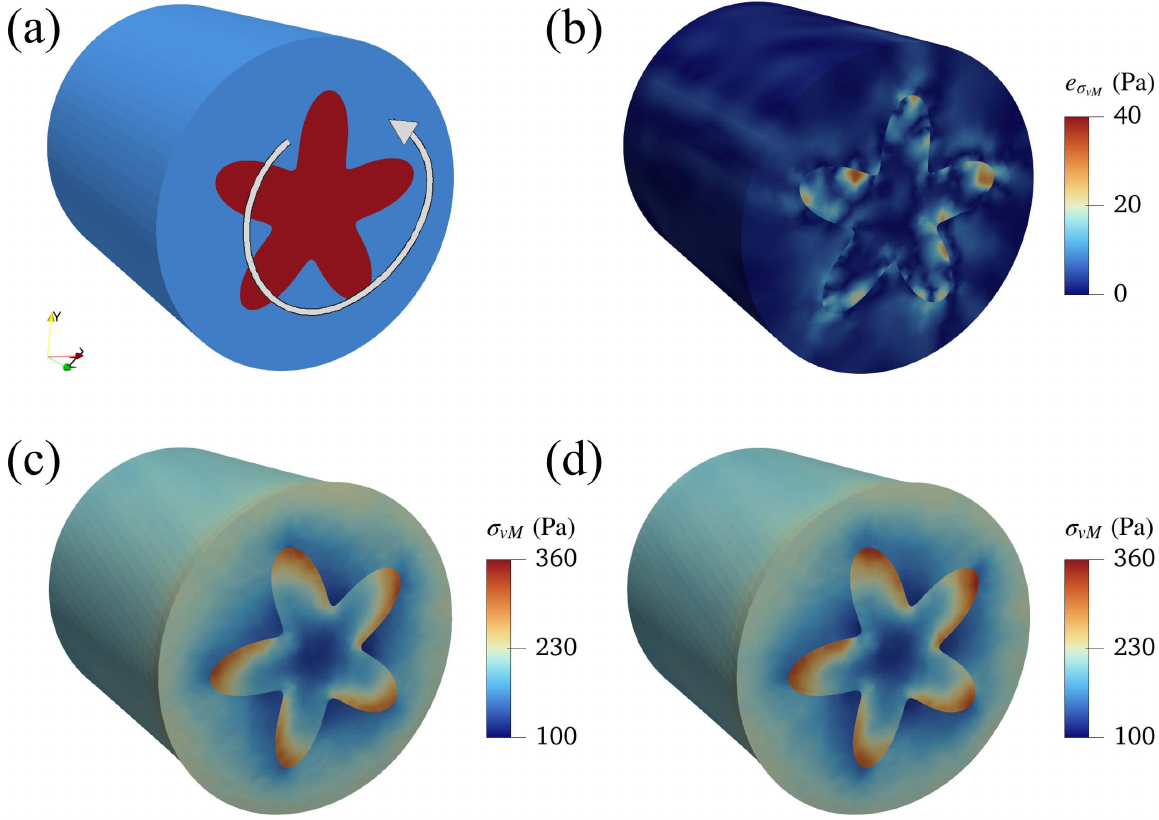}
    \caption{Hyperelastic composite rod with a five-lobed re-entrant inclusion under torsion. (a) Reference configuration showing the stiff inclusion (dark red) embedded in the compliant matrix (blue), with torque applied about the rod axis at \(z=H\) and the face \(z=0\) clamped. (b) Absolute error in the von Mises stress, \(e_{\sigma_\text{vM}}=\lvert\sigma_\text{vM}^{\mathrm{PI\text{-}GNN}}-\sigma_\text{vM}^{\mathrm{FEM}}\rvert\), between the \pignn prediction and the FE reference. (c) Deformed FE configuration and (d) corresponding \pignn prediction, both showing \(\sigma_\text{vM}\). Panels (c) and (d) use the same stress scale, while panel (b) uses an independent error scale. Stresses are reported in \si{\pascal}, and the deformed configurations are shown at a 1:1 scale.}
    \label{fig:rod_disp}
\end{figure}

\noindent The final example combines a re-entrant material interface and a fully 3D torsional stress state, on an unstructured mesh with a curved, traction-free lateral boundary. The geometry is a circular rod of radius \(R=\SI{1}{\milli\meter}\) and height \(H=\SI{2}{\milli\meter}\) (\(H/2R=1\)), containing a concentric five-lobed re-entrant inclusion extending over the full height of the rod. The phases are modeled as compressible Neo-Hookean solids using the same material properties as the hydrogel example. Reusing the material pair isolates the effects of 3D torsion and curved geometry while keeping the constitutive contrast fixed. The bottom face \(z=0\) is fully clamped, \(\mathbf{u}=\mathbf{0}\), while the top face \(z=H\) is subjected to a torque \(M_z=\SI{0.25}{\micro\newton\meter}\) (see \Cref{fig:rod_disp}a). The torque is chosen to produce finite deformation, with a mean twist exceeding \(\SI{20}{\degree}\), and the lateral surface is traction-free.

The conforming mesh contains \(38430\) nodes and \(214260\) tetrahedra, generated by extruding an adaptively refined cross-sectional triangulation containing \(1830\) nodes and \(3571\) triangles through 20 layers. The network architecture is identical to that used for the cube under torsion, with \(\vb{f}_i\in\mathbb{R}^5\) and \(\tilde{\vb{u}}_i\in\mathbb{R}^3\). The Dirichlet mask is applied component-wise on the clamped face, and the total energy is assembled element-wise over the tetrahedra, with the applied traction integrated over the triangular facets of the loaded face.

This example combines two aspects that have not previously been tested simultaneously. First, the unstructured graph is substantially denser than the cube's, with \(4.1\) times as many nodes and \(4.5\) times as many elements, providing a test of message passing on a locally refined 3D graph. Second, the applied load places the response well beyond the small-strain regime: the loaded face rotates by \(\SI{26.05}{\degree}\), while the rim shear is \(\gamma=\phi R/H=0.227\), approximately an order of magnitude larger than the \(0.02\)-\(0.05\) range for which the linear torsion solution remains applicable. Consequently, the observed warping and second-order axial stresses represent genuine finite-strain effects rather than corrections within the small-strain approximation.

\Cref{fig:rod_disp}c \& d compare the reference and predicted von Mises stress fields on the deformed configuration. The stiff inclusion attracts load, with \(\sigma_{\mathrm{vM}}\) reaching approximately \(\SI{360}{\pascal}\) along the lobe flanks, while the inclusion core remains at approximately \(\SI{100}{\pascal}\) owing to its proximity to the twist axis and correspondingly low shear. The matrix between the lobe tips and the free surface sustains intermediate stresses of approximately \(\SI{230}{\pascal}\). The \pignn captures the overall stress distribution, including the lobe peaks and steep gradients near the concave notches, with slightly smoother gradients in the surrounding matrix. The absolute error remains below \(\SI{10}{\pascal}\) at \(85\%\) of the nodes, corresponding to approximately \(3\%\) of the peak stress, and reaches a maximum of \(\SI{37.9}{\pascal}\) in localized regions near the lobe tips and re-entrant notches (see \Cref{fig:rod_disp}b). The global deformation is similarly well reproduced: the mean twist of the loaded face is \(\SI{25.876}{\degree}\), compared with \(\SI{26.051}{\degree}\) for the reference, corresponding to an error of \(0.7\%\), while the peak displacement magnitude is \(\SI{0.4714}{\milli\meter}\), compared with \(\SI{0.4745}{\milli\meter}\).

\Cref{tab:rod_error_stats} summarizes the error statistics. The displacement magnitude is recovered with \(L^{2}_{\lVert\vb{u}\rVert}=0.83\%\), while the in-plane components have errors of \(0.75\%\) and \(0.93\%\), with \(R^2\geq0.999\) for both. Thus, the curved free surface and re-entrant interface have little effect on the predicted twist response. The load-carrying shear component \(\sigma_{\theta z}\) has \(L^{2}_{\alpha}=5.33\%\), \(R^2=0.978\), and a median signed error of \(-0.41\%\). The von Mises stress shows similar accuracy, with \(L^{2}_{\alpha}=5.00\%\) and \(R^2=0.980\). Restricting the evaluation to the interface region yields essentially the same \(R^2\) values as the global evaluation, indicating that the near-interface region does not substantially reduce correlation.
\begin{table}[H]
    \centering
    \small
    \setlength{\tabcolsep}{14pt}
    \caption{Summary of error statistics for the composite rod under torsion. The \(L^{2}_\alpha\) values are volume-weighted, except for \(u_x\) and \(u_y\), which are evaluated nodally; the volume-weighted \(L^2\) error of the displacement vector is \(0.8305\%\).}
    \label{tab:rod_error_stats}

    \begin{tabular}{l S[table-format=2.4] S[table-format=-1.4] S[table-format=1.4] S[table-format=1.4]}
        \toprule \toprule
        & & & \multicolumn{2}{c}{\(R^2\)} \\
        \cmidrule(l){4-5}
        Field & {\(L^{2}_\alpha\) (\%)} & {Median \(\text{err}_{\alpha}\) (\%)} & {Global} & {Near-field} \\
        \midrule
        \(u_x\)                &  0.7482 & -0.4203 & 0.9999 & 0.9999 \\
        \(u_y\)                &  0.9258 &  1.1214 & 0.9999 & 0.9999 \\
        \(\sigma_{\theta z}\)  &  5.3262 & -0.4147 & 0.9777 & 0.9777 \\
        \(\sigma_\mathrm{vM}\)        &  5.0012 & -0.5091 & 0.9804 & 0.9805 \\
        \bottomrule \bottomrule
    \end{tabular}
\end{table}

\aashay{The energy balance provides an independent check unavailable from pointwise metrics. The converged \pignn attains a total potential energy of $\Pi = \SI{-5.978e-2}{\kilo\pascal\cubic\milli\meter}$ against $\SI{-6.022e-2}{\kilo\pascal\cubic\milli\meter}$ for the FE reference field evaluated with the same discrete functional, a relative difference of $+0.72\%$ within just 30,000 epochs. The prediction, therefore, approaches the reference minimum from above, as required by the principle of minimum potential energy. The model uses the same architecture and learning rate as the preceding examples. It takes \SI{37237}{\second} to train for 30,000 epochs on a single thread of \texttt{Intel Xeon(R) Gold 5220R}. Subsequently, returns the full field in a single forward pass in \SI{12.42}{\second}, against \SI{672.35}{\second} for the \texttt{FEniCSx} Newton\textendash Raphson solve. While the training time for PIGNN is $\approx55$ times FEM solve, the inference time is $\approx0.02$ times FEM solve, which shows this method is a promising candidate for parametric surrogate models where multiple inferences need to be made from a single trained model.}

\section{Summary and Outlook}\label{sec:conclusion}
\noindent \amiya{We presented a variational, label-free physics-informed graph neural network (\pignn) for two-phase heterogeneous solids under small-strain linear elasticity and finite-strain Neo-Hookean hyperelasticity in 2D and 3D. The solver minimizes the discrete total potential energy on a conforming mesh graph as a single scalar objective, with material heterogeneity carried at the element level. The interface emerges from the stiffness contrast between neighboring elements exactly as in FEM, with no penalty term and no prescribed transition width. The discrete energy on $\mathcal{P}_1$ elements coincides with the FE Ritz functional, and the converged FE field is the exact minimizer of the training loss. \aashay{Because $\nabla{u}$ is recovered from the $\mathcal{P}_1$ element shape functions rather than by autograd, the trial field carries no $C^1$ requirement\textemdash the constraint that obliges coordinate-based energy methods to adopt $\tanh$ or sinusoidal activations\textemdash a ReLU network is admissible here.}}

\amiya{Three sets of results establish the practical content of this construction. First, sweeping the stiffness contrast over $\Econ\in[10^{-2},10^{2}]$ with all settings frozen, the \pignn von Mises error remains below 3.58\% while a strong-form PINN with fixed transition width reaches 5.58\%; the \pignn error depends on the magnitude of the contrast but not its sign, whereas the PINN accuracy depends on both. Second, the \pignn roughly halves the $\sigma_{xx}$ error of an energy-based PINN on the perforated-plate benchmarks under discretization transfer \citep{Li2021}. \aashay{Third, ablating message passing while preserving the Ritz objective isolates what the graph contributes: at a fixed reduced training budget, the full model attains a 19-43\% (increasing with element size) lower error across a mesh sweep, making it a faster-converging model. The ablated version of our model reaches the minima of the loss function only with sufficient training budget, as discussed in (\ref{app:mp_converged}). Message passing, therefore, accelerates the approach to the shared Ritz minimizer rather than lowering it, and the advantage that persists at convergence is especially for the near-interface displacement field.}}

\aashay{The principal limitation is that the framework remains a per-instance solver: each new BVP requires separate training. The accuracy is bounded by the $\mathcal{P}_1$ Ritz minimizer, and the message-passing advantage narrows under strong mesh refinement and lenient training budgets. As shown in \Cref{sec:rod_torsion}, the training time exceeds a single FE solve by $55$ times. However, when we compare only the inference times, the PIGNN delivers the full field in a single forward pass, $50$ times faster than a \texttt{FEniCSx} Newton-Raphson solve. These times are case-specific, and the speed-up of inference depends mainly on the number of elements in the graph. Another systematic feature recurs throughout the suite: stresses are recovered with lower accuracy than displacements because they are determined \emph{post hoc} via element-wise constitutive evaluation.}

\amiya{Future work will upscale the solver to a parametric surrogate trained on families of geometry, loading, and material variations~\citep{Rezaei2022,Maurizi2022}, a setting in which the unsupervised energy objective eliminates the FE dataset that data-driven surrogates require. The graph architecture also extends to elastoplasticity~\citep{Niu2023} and anisotropic hyperelasticity~\citep{Abueidda2021}, and the self-supervised loss is well-suited to inverse problems such as parameter identification and elasticity imaging~\citep{Kamali2023}.}

\nolinenumbers
\noindent
\section*{Data availability}
\noindent The data and source code required to reproduce the results presented in this study will be openly available at \url{https://github.com/aashay-y/Variational-PIGNN-composites}. The repository contains the implementation, input files, and documentation necessary to reproduce the numerical experiments and figures reported in this work. 

\section*{Conflict of Interest}
\noindent The authors declare no conflicts of interest.

\appendix
\setcounter{figure}{0}
\section{Mesh Convergence study for Kirsch verification}\label{app:mesh_conv_1_hole}
\noindent Since the peak boundary stress from a \(\mathcal{P}_1\) (constant-strain) triangulation is mesh-dependent, we performed a native FEM mesh-convergence study before adopting the FEM fields as the verification benchmark. The benchmark is the classical single-hole problem: a $20\times20$ plane stress plate ($E=\SI{210}{\giga\pascal}$, $\nu=0.3$) with one central circular hole of radius $r=1$, subjected to uniform traction $\sigma_0=\SI{1}{\mega\pascal}$ on the right edge (see~\Cref{sec:single_hole}). 

\begin{figure}[ht]
    \centering
    \includegraphics[width=0.9\linewidth]{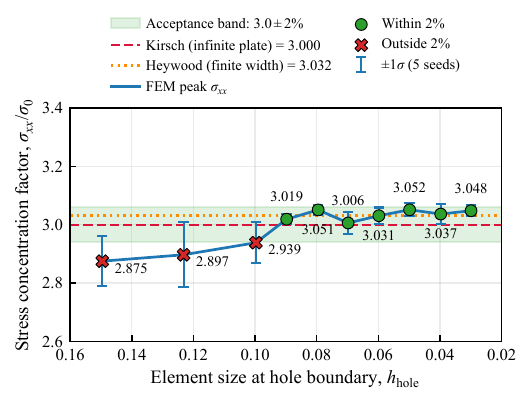}
    \caption{Mesh convergence of the SCF for the single-hole verification example (plane stress linear elasticity; $r=1$, plate $=20$, $\nu=0.3$, gross traction $\sigma_0=\SI{1}{\mega\pascal}$). Markers show the seed-averaged FEM peak nodal $\sigma_{xx}/\sigma_0$ \textit{vs.} the element size at the hole boundary $h_{\text{hole}}$ (coarse \(\rightarrow\) fine). The error bars denote the first standard deviation over five independent mesh seeds. The shaded region is the $3.0\pm2\%$ acceptance band $[2.94,3.06]$; green circles lie inside the band (acceptable), and red crosses fall outside it. Horizontal lines denote the Kirsch infinite-plate solution ($3.0$, dashed) and the Heywood finite-width value ($3.032$, dotted). The SCF converges into the band for $h_{\text{hole}}\leq0.09$, plateauing at $\approx3.04$.}
    \label{fig: mesh_convergence_single_hole}
\end{figure}

For a linear elastic problem, the stress concentration factor (SCF) is independent of $E$ and the load magnitude; the peak stress $\sigma_{xx}$ equals the SCF. We refined the characteristic element size at the hole boundary, $h_{\text{hole}}=2\pi r/(\text{number of hole-boundary facets})$, from $0.03$ to $0.15$, while holding the far-field grading fixed, such that only the near-hole resolution is varied. The SCF is taken as the seed-averaged nodal peak $\sigma_{xx}/\sigma_0$ value. Here we used the same area-weighted cell$\rightarrow$nodal projection as used in the analytical/FEM/\pignn comparison pipeline. As the unstructured Poisson-disk places interior nodes randomly, each element size is meshed with five independent random number generator seeds and reported as mean $\pm$ one standard deviation. The per-size scatter collapses from $\sim0.07$–$0.11$ ($h_{\text{hole}}\geq0.10$) to $0.018$–$0.039$ ($h_{\text{hole}}\leq0.09$), confirming the fine meshes are stable. In~\Cref{fig: mesh_convergence_single_hole}, the SCF increases from $2.88$ (at $h_{\text{hole}}=0.15$) to plateau at $\approx3.04$, consistent with the Kirsch infinite-plate value of $3.000$ and the Heywood finite-width estimate of $3.032$ ($d/W=0.1$). Every mesh with $h_{\text{hole}}\leq0.09$ falls inside the acceptance band ($3.0\pm2\%$), while coarser meshes are under-resolved and fall low. 

The independent Kirsch check is also satisfied $\sigma_{yy}\rightarrow-1$ at the hole side, and the peak $\sigma_{xx}$ is located on the hole boundary at the crown. For the \pignn training mesh, with \(h_{\mathrm{hole}}\approx 0.025 \) and 3577 nodes, the computed SCF is 3.04, which lies within the prescribed acceptance band. This confirms that the FE reference solution used \aashay{for verification} is sufficiently mesh-converged with respect to the SCF.

\section{Architecture Ablation Study}\label{app:ablation}
\noindent An architecture ablation study is conducted to quantify the trade-off between prediction accuracy and model complexity for the proposed \pignn model. The ablation study is performed on a stiff re-entrant inclusion, while keeping all other training settings fixed, including the learning rate, number of training epochs, material properties, dataset, and random seed (see~\Cref{sec:irregular-inclusion}). Only the network architecture is varied by changing the number of message-passing layers and the width of the hidden neurons.

\begin{figure}[ht]
    \centering
    \includegraphics[width=0.75\linewidth]{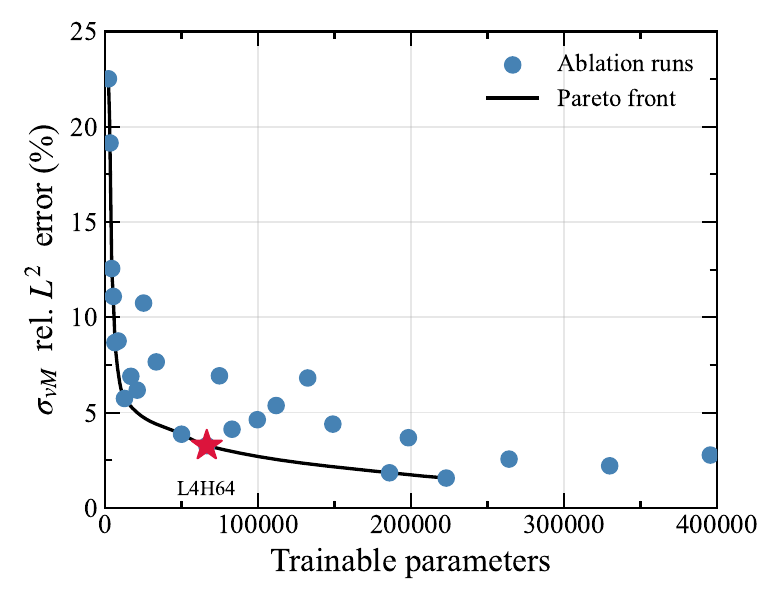}
    \caption{Trade-off between model complexity (trainable parameters) and area-weighted \(L^{2}\) for the \pignn architecture sweep. Blue markers denote the evaluated configurations, the solid curve represents the Pareto front, and the selected architecture (\texttt{L4H64}) is labeled with a red star.}
    \label{fig:ablation_pareto}
\end{figure}

Each architecture is evaluated using two competing objectives: (i) the area-weighted \(L^{2}\) of the von Mises stress field, which measures prediction accuracy, and (ii) the number of trainable parameters, which serves as a measure of model complexity.~\Cref{fig:ablation_pareto} shows the resulting architectures and the corresponding Pareto front. The Pareto front identifies architectures for which no further reduction in prediction error is possible without increasing model complexity. Based on this trade-off, \aashay{the configuration \texttt{L4H64} is selected for all application examples}, as it provides a favorable balance between accuracy and computational cost.

\section{Stiffness-contrast sweep}\label{app:sweep}
\noindent \amiya{The sweep in~\Cref{sec:pinn_vs_pignn_inclusion} considers twenty-five logarithmically spaced stiffness contrasts spanning $\Econ \in [10^{-2}, 10^{2}]$ (with twelve contrasts on either side of the unit-contrast homogeneous limit), obtained by varying \(E_{\mathrm{inc}}\) while keeping \(E_{\mathrm{mat}}=\SI{1500}{\mega\pascal}\), \(\nu_{\mathrm{mat}}=\nu_{\mathrm{inc}}=0.4\), the geometry, and the applied traction \(T=\SI{0.025}{\mega\pascal}\) fixed. The PINN uses four hidden layers of 64 neurons with \(\tanh\) activations, a \(128\times128\) collocation grid, a transition width \(\delta=0.01\), and \(5000\) \texttt{BFGS} iterations. The \pignn uses the \texttt{L4H64} configuration, trained for \(50{,}000\) \texttt{Adam} epochs at \(\eta=10^{-3}\) (see~\ref{app:ablation}). Both solvers are trained from three independent random initializations at each contrast. The sweep, therefore, compares the two methods under fixed configurations rather than isolating the effect of any individual hyperparameter.}

Both solvers are evaluated using the same implementation on a common held-out mesh with \(3034\) nodes and \(5944\) elements, which is also used for the FE reference solution. Although the PINN is trained without a mesh, its coordinate network is evaluated at the nodes of this mesh, allowing the element areas to be used consistently for both methods in~\Cref{eq:4.2}. This area weighting is important because the element areas vary by a factor of \(168\), with the smallest elements concentrated near the material interface, where the errors are also largest. For the \pignn at a contrast of \(100\), removing the area weights increases the reported relative \(L^2\) error from \(3.58\%\) to \(7.96\%\). The weighting, therefore, materially affects the reported error, making its consistent use for both methods essential for a meaningful comparison.

\subsection{The unit contrast limit}\label{app:ratio_one}
\noindent For unit contrast \(E_{\mathrm{inc}}/E_{\mathrm{mat}}=1\), the PINN error is \(0.002\%\), compared with \(0.598\%\) for the \pignn. This exceptionally small PINN error is specific to the homogeneous case and follows from its trial-function construction: the prescribed background field already coincides with the exact uniform stress solution, leaving only a vanishing correction for the network to learn. The \pignn has no analogous prescribed stress field and must recover the equilibrium state by minimizing the discrete potential energy. Thus, this isolated result does not indicate a general advantage in accuracy for the PINN and is not representative of heterogeneous cases.

The residual \pignn error is primarily due to incomplete optimization rather than to spatial discretization. For all three seeds, the stored network solutions have higher discrete energy than the FE solution, confirming that the discrete minimizer has not been reached. The corresponding element-wise error fields are essentially uncorrelated (correlation coefficient = \(0.013\), \(0.013\), and \(0.026\)), further indicating optimization noise rather than a systematic spatial error. Consistent with this interpretation, reducing the learning rate to \(10^{-4}\) and selecting the lowest-energy checkpoint reduces the error from \(0.598\%\) to \(0.368\%\).

The residual error also exhibits small-scale roughness in the predicted displacement field. A local least-squares plane fit, which preserves linear fields exactly, reduces the von Mises error by \(33\)--\(44\%\) after one pass and by \(71\%\) after eight passes, while leaving the FE solution unchanged. This confirms that part of the residual originates from optimization-induced displacement roughness. Such smoothing is appropriate only for this homogeneous case, where the exact solution is linear, and it would not preserve the genuine gradients and stress discontinuities that arise at material interfaces. Therefore, the exceptionally small PINN error at unit contrast and its increasing error with stiffness contrast arise from distinct mechanisms rather than a single monotonic difference in solver accuracy.

\subsection{Comparison of converged solutions}\label{app:mp_converged}
\noindent \aashay{The \Cref{sec:message_passing} uses a fixed reduced budget of 5000 epochs since the objective was to compare the convergence speed and not the converged solutions. Since the models minimize the same functional, they can be compared directly once the optimization has converged. Here the same reduced model setting (\texttt{L4H32}, 17,250 parameters, \(\eta=10^{-3}\)) but trained to 30,000 \texttt{Adam} epochs. The mesh is similar; fixed interface-conforming mesh (3034 nodes, 5944 elements) and for a fixed contrast \(\Econ = 10\).}
\begin{figure}[H]
    \centering
    \includegraphics[width=0.6\linewidth]{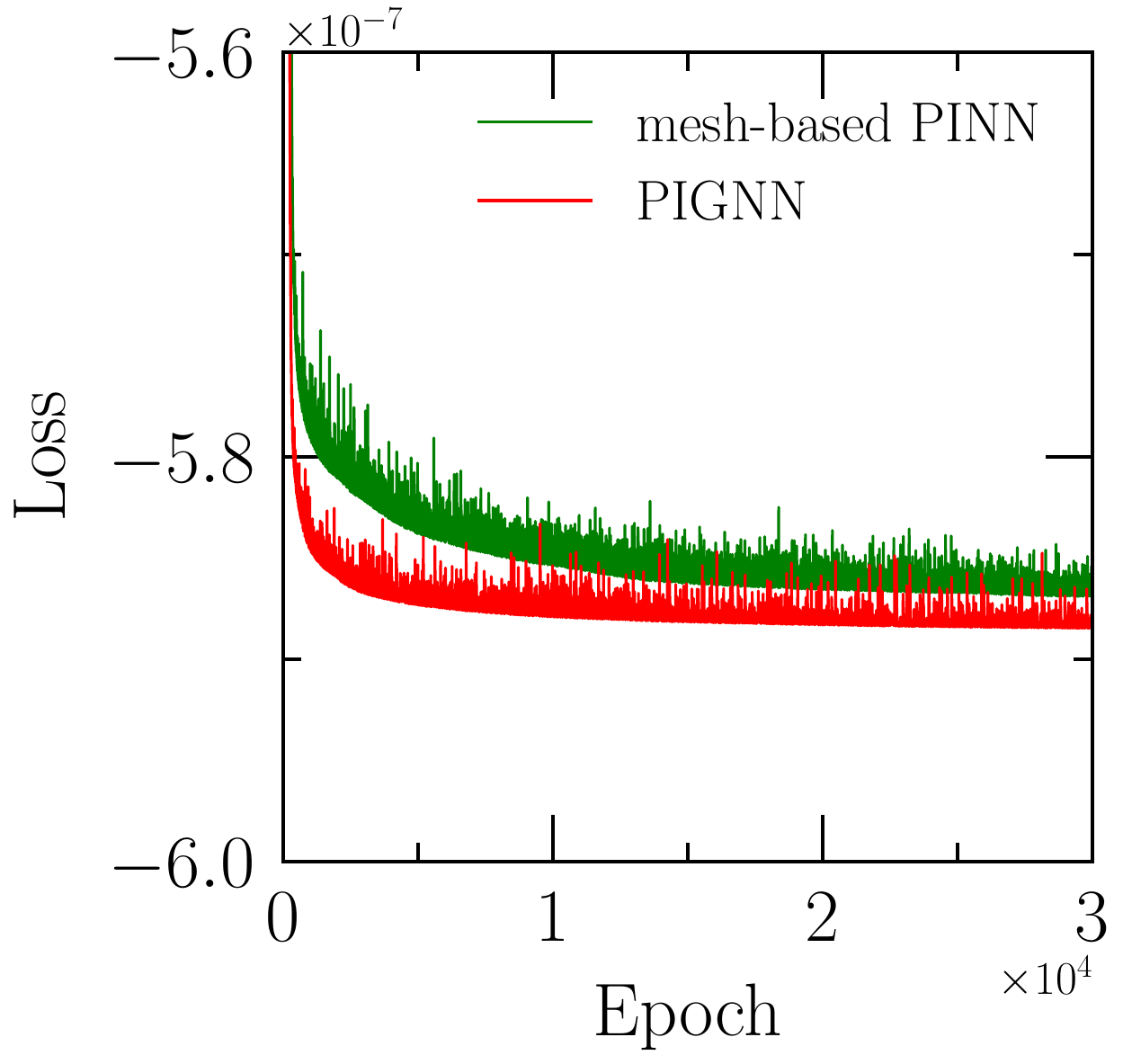}
    \caption{the total potential \(\Pi\) descends faster for \pignn in the first \(\sim\)3000 epochs. It reaches a \(\Pi\) value that the mesh-based PINN only attains near 30,000 epochs. The ablated model curve then slowly closes the gap.}
    \label{fig:loss_30k}
\end{figure}

\aashay{The converged energies are \(\Pi = -5.885\times10^{-7}\) (\pignn, epoch 29,997) and \(\Pi = -5.869\times10^{-7}\) (mesh-based PINN, epoch 29,818). The FE reference \(\Pi_{\mathrm{FEM}} = -5.889\times10^{-7}\) \si{\mega\pascal\milli\meter\squared}. Therefore, we can consider that the models reach the minima after sufficient epochs, but PIGNN converges faster.}

\aashay{Comparing error metrics for converged fields we get, \(L^{2}_{\sigma_\text{vM}} = 4.60\%\) versus \(4.64\%\) so the advantage seen at 5000 epochs has closed. Displacement retains a gap of \(L^{2}_{\mathrm{u}} = 0.239\%\) versus \(0.709\%\) while the interface band \(0.250\%\) versus \(0.837\%\) and near-field \(R^{2}_{u_x} = 0.9997\) versus \(0.9985\).}

\bibliographystyle{elsarticle-num-names}
\biboptions{comma,round}
\bibliography{Ref}
\end{document}